\documentstyle[11pt]{article}
 
\newcommand{\ncm}{\newcommand}
\ncm{\rncm}{\renewcommand}
 
\ncm{\lb}[1]{\label{#1}}
\rncm{\sec}{\setc{0}\section}
\ncm{\beq}{\begin{equation}}
\ncm{\eeq}{\end{equation}}
\ncm{\bea}{\begin{eqnarray}}
\ncm{\beanon}{\begin{eqnarray*}}
\ncm{\eea}{\end{eqnarray}}
\ncm{\eeanon}{\end{eqnarray*}}
 
\rncm{\theequation}{\thesection.\arabic{equation}}
\ncm{\setc}[1]{\setcounter{equation}{#1}}
\newcounter{eqnr}
 
\newenvironment{eqnarrayabc}{\stepcounter{equation}
  \setcounter{eqnr}{\value{equation}}\setc{0}
  \rncm{\theequation}{\thesection.\arabic{eqnr}\alph{equation}}
  \begin{eqnarray}}{\end{eqnarray}\setc{\value{eqnr}}}
 
\ncm{\eql}{&\quad\stepcounter{equation}(\theequation a)
           &\qquad\qquad}
\ncm{\eqr}{&\quad(\theequation b)}
 
\newenvironment{lrarray}{$$\begin{array}{llcrr}}
                        {\end{array}$$}
 
\ncm{\beabc}{\begin{eqnarrayabc}}
\ncm{\eeabc}{\end{eqnarrayabc}}
\ncm{\eqboxabc}[3]{\newline\parbox[t]{1.5cm}{#1}\hfill
  \parbox[b]{12cm}{\begin{eqnarray*} #3\end{eqnarray*}}\hfill
   \parbox[b]{1.5cm}{\vspace{-0.0cm}
  \begin{eqnarrayabc}#2\end{eqnarrayabc}}\newline}
\ncm{\nn}{\nonumber \\}
 
\newtheorem{thm}{Theorem}[section]
\newtheorem{prop}[thm]{Proposition}
\newtheorem{lem}[thm]{Lemma}
\newtheorem{scho}[thm]{Scholium}
\newtheorem{defi}[thm]{Definition}
\newtheorem{coro}[thm]{Corollary}
\newtheorem{exam}[thm]{Example}
 
\def\End{\mbox{End}\,}
\def\Ker{\mbox{Ker}\,}
\rncm{\Im}{\mbox{Im}}
\def\Span{\mbox{\rm Span}}
\def\Aut{\mbox{Aut}\,}
\def\Out{\mbox{Out}\,}
\def\Hom{\mbox{Hom}\,}

\def\id{\mbox{id}\,}
\def\Mat{\mbox{Mat}\,}
\def\Center{\mbox{Center}\,}
\def\Hypercenter{\mbox{Hypercenter}\,}
\def\Ad{\mbox{Ad}\,}

\def\Inv{\mbox{\rm Inv}\,}
\def\Coinv{\mbox{\rm Coinv}\,}
\def\A{{\cal A}}

\def\CA{\mbox{Center}\,A}

\def\duZ{{\hat Z}}

\def\D{{\cal D}}

\def\I{{\cal I}}
\def\LI{{\cal I}^L}
\def\RI{{\cal I}^R}
\def\duLI{{\hat\LI}}

\def\IL{{\cal I}^L(A)}
\def\IR{{\cal I}^R(A)}
\def\duIL{{\cal I}^L(\duA)}

\def\Ind{\mbox{Index}\,}

\def\1{1}
\def\du1{{\hat\1}}
\def\o{\otimes}
\def\x{\times}
\def\cop{\Delta}
\def\ducop{{\hat\cop}}
\def\eps{\varepsilon}
\def\dueps{{\hat\eps}}
\def\duS{{\hat S}}
\def\duA{{\hat A}}
\def\tr{\mbox{tr}\,}
\def\Tr{\mbox{Tr}\,}
\def\c{_{(1)}}
\def\cc{_{(2)}}
\def\ccc{_{(3)}}

\def\cp{_{(1')}}
\def\ccp{_{(2')}}
\def\PL{\sqcap^L}
\def\PR{\sqcap^R}
\def\PLR{\sqcap^{L/R}}
\def\duPL{{\hat\PL}}
\def\duPR{{\hat\PR}}

\def\Lra{\Leftrightarrow}
 
\def\cros{\raise1.9pt\hbox{$\scriptscriptstyle
          > $}\!\raise1.5pt\hbox{$\scriptstyle\triangleleft\,$}}

\def\C{\,{\raise1.5pt\hbox{$\scriptscriptstyle |$}
        \thinmuskip=4mu \!\!C\thinmuskip=3mu}}
\def\Z{{Z\!\!\!Z}}
\def\N{{\thinmuskip = 5.5mu I\!N\thinmuskip = 3mu}}

\def\bra{\langle}
\def\ket{\rangle}
\def\la{\!\rightharpoonup\!}
\def\ra{\!\leftharpoonup\!}
\def\qed{\hfill {\it Q.e.d.}}
\def\Proof{{\em Proof}\,:\ }
\def\triv{D_{\varepsilon}}
\def\V{V_{\eps}}
\ncm{\Sec}{{\cal S}\!{\it ec}\,}
\ncm{\Vac}{{\cal V}\!{\it ac}\,}
\ncm{\Hyp}{{\cal H}\!{\it yp}\,}

\ncm\amalgo[1]{{\lower9pt\hbox{$\o$}\atop\raise2pt\hbox{
            $\scriptscriptstyle #1$}}}
 
\def\amo{{\lower9pt\hbox{$\otimes$}\atop\raise2pt\hbox{
          $\scriptscriptstyle A^R\equiv\hat A^L$}}}

\begin{document}
 
\large
\title{\bf Weak Hopf Algebras \\
       I. Integral Theory and $C^*$-structure}
 
\author{\sc Gabriella B\"ohm $^1$,  Florian Nill $^2$,
            Korn\'el Szlach\'anyi $^3$ \\}

\date{}
 
\maketitle
 
\footnotetext[1]{
Research Institute for Particle and Nuclear Physics, Budapest,
H-1525 Budapest 114, P.O.B. 49, Hungary\\
E-mail: BGABR@rmki.kfki.hu\\
  Supported by the Hungarian Scientific Research Fund, OTKA --
  T 016 233}
 
\footnotetext[2]{
Institut f\"ur Theoretische Physik, FU--Berlin
Arnimallee 14, D-14195 Berlin, Germany\\
E-mail: NILL@mail.physik.fu-berlin.de
\\Supported by the DFG, SFB 288 "Differentialgeometrie und
Quantenphysik"}
 
\footnotetext[3]{
Research Institute for Particle and Nuclear Physics, Budapest,
H-1525 Budapest 114, P.O.B. 49, Hungary\\
E-mail: SZLACH@rmki.kfki.hu\\
  Supported by the Hungarian Scientific Research Fund,
  OTKA -- T 020 285.}
 
\vskip 2truecm

\begin{abstract}
We give an introduction to the theory of
weak Hopf algebras
proposed recently as a {\em coassociative} alternative of weak
quasi-Hopf algebras. We follow an axiomatic approach keeping as
close as possible to the "classical" theory of Hopf algebras.
The emphasis is put on the new structure related to the presence
of
canonical subalgebras $A^L$ and $A^R$ in any weak Hopf algebra $A$
that play the role of non-commutative numbers in many respects.
A theory of integrals is developed in which we show how the
algebraic properties of $A$, such as the Frobenius property, or
semisimplicity, or innerness of the square of the antipode, are
related to the existence of non-degenerate, normalized, or Haar
integrals. In case of
$C^*$-weak Hopf algebras we prove the existence of a unique Haar
measure $h\in A$ and of a canonical grouplike element
$g\in A$ implementing the square of the antipode and factorizing
into left and right elements $g=g_L\,g_R^{-1}$, $g_L\in A^L$,
$g_R\in A^R$. Further discussion of the $C^*$-case will be
presented in Part II.
\vskip 1.2truecm\small
\hskip 8truecm {\em To appear in J. Algebra} 
\end{abstract}

\vfill\eject
\small
\tableofcontents
\normalsize
\section{Introduction}

Weak Hopf algebras have been proposed recently \cite{BSz,Sz,WBA}
as a new
generalization of ordinary Hopf algebras that replaces Ocneanu's
paragroup \cite{Ocneanu 1}, in the depth 2 case, with a
concrete "Hopf algebraic" object. The earlier proposals of {\em
face algebras} \cite{Hayashi} or {\em quantum groupoids} \cite{Ocneanu 2}
are actually weak Hopf algebras even if not the most general
ones.
Also, the (finite dimensional) {\em generalized Kac algebras} of
T. Yamanouchi \cite{Yamanouchi} are weak Hopf algebras in our
sense \cite{WBA}, albeit with an involutive antipode.
 
In contrast to other Hopf algebraic
constructions such as the quasi Hopf algebras \cite{Drinfeld} or
the weak quasi Hopf algebras and rational Hopf algebras
\cite{MS,V,FGV} weak Hopf algebras are {\em coassociative}.
This allows one to define actions, coactions, and crossed
products as easily as in the Hopf algebra case. On the other
hand weak Hopf algebras have "weaker" axioms related to the unit
and counit: The comultiplication is non-unital, $\cop(\1)\neq
\1\o\1$ (like in weak quasi
Hopf algebras) and the counit is only "weakly" multiplicative,
$\eps(xy)=\eps(x\1\c)\eps(\1\cc y)$. This kind of "weakness" is
the "strength" of weak Hopf algebras because it allows (even in
the finite dimensional and semisimple case) the weak Hopf algebra
to possess non-integral (quantum) dimensions.
 
Thus weak Hopf algebras are not special cases of weak quasi
Hopf algebras and also not more general than them.
Nevertheless, in situations where only the representation category
of the quantum group matters, these two concepts are equivalent.
This is, of course, not surprising in view of MacLane's Theorem on
the equivalence of relaxed and strict monoidal categories
\cite{McLane}. In
fact not all of the potential of this theorem is utilized by weak
Hopf algebras because their representation category is not
quite strict: Only the associator is trivial but not the left and
right isomorphisms of the monoidal unit. Although a general
analysis clarifying the role of representation categories of
weak Hopf algebras within the set of monoidal categories is still
missing the examples constructed in \cite{BSz} using Ocneanu's
cocycle suggest that they
play a rather fundamental role, as long as they can accomodate to
arbitrary 6j-symbols.
 
So far weak Hopf algebras have been considered only under the
additional assumption of finite dimensionality. Although a good
deal of the results can be generalized to the infinite dimensional
case, finite dimension is particularly attractive because it
implies selfduality. Just like finite Abelian
groups or finite dimensional Hopf algebras, the finite
dimensional weak Hopf algebras (WHA) are {\em
selfdual} in the following sense.
If $A$ is a WHA then its dual space $\duA$ is canonically
equipped with a weak Hopf algebra structure. Furthermore this
duality is reflexive, $(\duA)\hat{\mbox{~}} \cong A$. This is a
feature which makes WHAs more natural objects of study than either
finite (non-Abelian) groups or finite dimensional (weak) quasi
Hopf algebras.
 
The main motivation for studying WHAs comes from quantum field
theory and operator algebras and consists roughly of the following
two symmetry problems.
\begin{description}
\item[I.]
If $N\subset M$ is an inclusion of algebras
satisfying
certain conditions then find a (unique) "quantum group" $G$ and an
action of $G$ on $M$ such that $N=M^G$, the invariant subalgebra.
\item[II.]
The dual problem is to find a "quantum group" $\hat G$
acting on $N$ such that $M$ is isomorphic to the crossed product
$N\cros\hat G$.
\end{description}
Of course, determining the appropriate notion of "quantum group",
as well as its action, is part of the problem. If $N\subset M$ is
a finite index irreducible depth 2 inclusion of von Neumann
factors then the answer is known by \cite{Longo} to be a finite
dimensional $C^*$-Hopf algebra. In \cite{NSzW} we will show that if we
allow the inclusion to be reducible and $N$ and $M$ to have
arbitrary finite dimensional centers then the appropriate
"quantum group" is a $C^*$-weak Hopf algebra.
Even in case of inclusions of certain associative (non-$^*$)
algebras
the notion of a WHA over an arbitrary field $K$, introduced in
this paper, may provide a useful invariant.
 
In Section 2 we introduce the axioms of weak bialgebras and
weak Hopf algebras over a field $K$
and discuss their consequences. If $K=\C$, the complex field, then
these axioms are equivalent to those of \cite{Sz}. The present axioms
have the advantage of being manifestly selfdual and almost each of
them having an ancestor among the Hopf algebra axioms which it
generalizes. In discussing the consequences particular attention
is paid to the canonical subalgebras $A^L$ and $A^R$ present in
any WHA both of which reducing to the scalars $K\1$ if $A$ is a
Hopf algebra. From many point of views these subalgebras behave
like non-commutative generalizations of numbers. Just to mention
some: 1. $A^L$ and $A^R$ are separable $K$-algebras. 2. The
trivial left $A$-module is a representation on the $K$-space $A^L$
(or on $A^R$). 3. The dual weak Hopf algebra $\duA$ have left and
right subalgebras $\duA^L$ and $\duA^R$ that are isomorphic to
$A^R$ and $A^L$, respectively. Of course, in order to realize the
idea of $A^L$ and $A^R$ being "non-commutative numbers" one should
completely get rid of the field $K$ from the outset. As yet we
have no concrete proposal for this scenario.

Section 3 is devoted to the study of integrals in weak Hopf
algebras. Using the notion of {\em weak Hopf modules} which is a
generalization of the Hopf modules \cite{Abe,Sweedler}
we show
that non-zero integrals exist. A weak Hopf version of Maschke's
Theorem characterizes semisimple WHAs as those possessing
normalized integrals. An other important class of WHAs are
those which are Frobenius algebras. They are characterized by
possessing non-degenerate left integrals. This class is a
selfdual class by the Duality Theorem of non-degenerate integrals.
We conclude with giving necessary and sufficient criteria for the
existence of Haar integrals, i.e. normalized non-degenerate
2-sided integrals in a WHA.
 
Section 5 contains the basic properties of weak $C^*$-Hopf
algebras such as the existence of a Haar
integral $h$ and a canonical grouplike element $g\geq 0$
implementing $S^2$ and the modular automorphism of the Haar
measure. As a consequence of the existence of Haar measures
the dual of a $C^*$-weak Hopf algebra is a $C^*$-weak Hopf algebra
again. Further analysis of $C^*$-WHAs will be given in Part
II where we discuss the representation category
and a notion of dimension which turns out to be non-commutative
in case of solitonic representations \cite{BNSzII} .

\sec{The Weak Hopf Calculus}
\subsection{The axioms}
 
\begin{defi}
A {\em weak bialgebra} (WBA) is a quintuple
$(A,\mu,u,\cop,\eps)$ satisfying Axioms 1, 2, and 3 below.
If $(A,\mu,u,\cop,\eps,S)$ satisfies Axioms 1, 2, 3, and 4 below
it is called a {\em weak Hopf algebra} (WHA).
\end{defi}
\begin{description}
\item [\bf Axiom 1.] $A$ is a finite dimensional associative
algebra over a field $K$ with multiplication
$\mu\colon A\o A\to A$ and unit $u\colon K\to A$. I.e. $\mu$ and
$u$ are $K$-linear and satisfy
\begin{description}
\item[\sf Associativity:] $\quad\mu\circ(\mu\o\id)=
                           \mu\circ(\id\o\mu)\quad$\hfill$(A.1)$
\item[\sf Unit property:] $\quad\mu\circ(u\o\id)=\id=
                           \mu\circ(\id\o u)\quad$\hfill$(A.2)$
\end{description}
(Later on we will suppress $\mu$ and $u$, just write $xy$ for
$\mu(x,y)$ and use the {\em unit element} $\1:=u(1)$ instead of
$u$.)
\item[\bf Axiom 2.] $A$ is a
coalgebra over $K$ with {\em comultiplication}
$\cop\colon A\to A\o A$ and {\em counit}
$\eps\colon A\to K$. I.e. $\cop$ and $\eps$ are $K$-linear and
satisfy
\begin{description}
\item[\sf Coassociativity:]
$\qquad(\Delta\o\id)\circ\Delta\ =\ (\id\o\Delta)\circ\Delta$
\hfill$(A.3)$
\item[\sf Counit property:]
$\qquad(\varepsilon\o\id)\circ\Delta\ =\ \id\
           =\ (\id\o\varepsilon)\circ\Delta$ \hfill$(A.4)$
\end{description}
\item[\bf Axiom 3.] For compatibility of the algebra and coalgebra
structures we assume
\begin{description}
\item[\sf Multiplicativity of the coproduct:] For all $x,y\in A$
 
\smallskip
$\Delta(xy)\ =\ \Delta(x)\Delta(y)$\hfill $(A.5)$
 
\smallskip
\item[\sf Weak multiplicativity of the counit:] For all $x,y,z\in
A$
 
\smallskip
$\eps(xyz)=\eps(xy\c)\eps(y\cc z)$\hfill$(A.6a)$
 
$\eps(xyz)=\eps(xy\cc)\eps(y\c z)$\hfill$(A.6b)$
 
\smallskip
\item[\sf Weak comultiplicativity of the unit:]\hfill
 
\smallskip
$\cop^2(\1)=(\cop(\1)\o\1)(\1\o\cop(1))$\hfill$(A.7a)$
 
$\cop^2(\1)=(\1\o\cop(1))(\cop(\1)\o\1)$\hfill$(A.7b)$
 
\smallskip
\end{description}
\item[Axiom 4.] There exists a $K$-linear map $S\colon A\to A$,
called the {\em antipode}, satisfying the following
\begin{description}
\item[\sf Antipode axioms:] For all $x\in A$
 
\smallskip
$x\c S(x\cc)=\eps(\1\c x)\1\cc$\hfill$(A.8a)$
 
$S(x\c)x\cc=\1\c\eps(x\1\cc)$\hfill$(A.8b)$
 
$S(x\c)x\cc S(x\ccc)=S(x)$\hfill$(A.9)$
 
\smallskip
\end{description}
\end{description}
In eqs. (A.6-9) we used a standard suffix notation for (iterated) coproducts,
omitting as usual summation indices and a summation symbol.
 
In the terminology of \cite{WBA} $(A,\mu,u,\Delta,\varepsilon)$ is called
a {\em weak bialgebra} if it satisfies the Axioms (A.1-5). There a weak
bialgebra is called {\em monoidal} if it satisfies (A.6) and it is called
{\em comonoidal} if it satisfies (A.7). As has been explored in detail in
\cite{WBA}, these (co)monoidality axioms are precisely designed to render
the category of $A$-modules (the category of $A$-comodules, respectively)
monoidal.
 
The {\em dual} of a weak bialgebra (weak Hopf algebra) $A$ is the
dual space $\duA:=\Hom_K(A,K)$ equipped with structure maps
$\hat\mu,\hat u,\ducop,\dueps\ (,\duS)$ defined by transposing
the structure maps of $A$ by means of the canonical pairing $\bra\
,\ \ket\colon\duA\times A\to K$ :
\beanon
\bra\varphi\psi,x\ket&:=&\bra\varphi\o\psi,\cop(x)\ket\\
\bra\du1,x\ket&:=&\eps(x)\\
\bra\ducop(\varphi),x\o y\ket&:=&\bra\varphi,xy\ket\\
\dueps(\varphi)&:=&\bra\varphi,\1\ket\\
\bra\duS(\varphi),x\ket&:=&\bra\varphi,S(x)\ket
\eeanon
where $\varphi,\psi\in\duA$ and $x,y\in A$.
 
Let $f$ and $g$ be maps from the $m$-fold tensor product $A^{\o
m}$ to the $n$-fold tensor product $A^{\o n}$ such that they are
composites of tensor products of the structure maps
$\mu,u,\cop,\eps,S$ and of the twist maps $\tau_{ij}$
interchanging the $i$-th and the $j$-th $A$ factors. Then the
equality $f=g$ is called an {\em $A$-statement}. Similarly one
defines the {\em $\duA$-statements}. Now every $A$-statement $Q\
::\ f=g$ determines an equivalent $\duA$-statement $Q^T\ ::\
f^T=g^T$ obtained by reversing the order of composition and
replacing $\mu$ with $\ducop$, $u$ with $\dueps$, $\cop$ with
$\hat\mu$, $\eps$ with $\hat u$, and $S$ with $\duS$. The
statement $Q^T$ is called the {\em transpose} of $Q$.If we now
substitute $\mu,u,\cop,\eps,S$, respectively in place of
$\hat\mu,\hat u,\ducop,\dueps,\duS$ in the statement $Q^T$ we
obtain a new $A$-statement $Q^\sim\ ::\ f^\sim=g^\sim$
which is not equivalent to $Q$ in general. This $Q^\sim$ will
be called the {\em dual} of $Q$. For example one can easily verify
that the WBA axioms satisfy $(A.1)^\sim=(A.3)$, $(A.2)^\sim
=(A.4)$, $(A.5)^\sim=(A.5)$, $(A.6a)^\sim=(A.7a)$, and
$(A.6b)^\sim=(A.7b)$. Thus the weak bialgebra
axioms form a selfdual set of statements. This implies that the
dual of a WBA is a WBA, too. The same holds for weak Hopf
algebras, since each one of the antipode axioms is a selfdual
statement. As a consequence of selfduality if $Q$ is a true
statement in a WBA or in a WHA then $Q^\sim$ is also true there. This
principle extends also to statements involving both $A$ and $\duA$
structure maps and canonical pairing(s).
 
As has been proven in \cite{WBA},
the above selfdual set of WHA axioms are equivalent
to the non-selfdual set of axioms given in \cite{Sz}.
In this work we will gradually reproduce all
axioms of \cite{Sz} as a consequence of the present ones.
 
For a weak Hopf algebra $(A,\1,\cop,\eps,S)$ the following
conditions are equivalent
\begin{itemize}
\item $A$ is a Hopf algebra;
\item $\cop(\1)=\1\o\1$\ ;
\item $\eps(xy)=\eps(x)\eps(y)$\ ;
\item $S(x\c)x\cc=\1\eps(x)$\ ;
\item $x\c S(x\cc)=\1\eps(x)$\ .
\end{itemize}
The proof of these assertions are either trivial or will become
trivial after acquainting the weak Hopf calculus developed in the
next subsections, see also \cite{WBA}.
 
\subsection{Weak bialgebras}
In a WBA define the maps $\PL,\PR\colon A\to A$ by the formulae
\beq                              \label{Pi}
\PL(x):=\eps(\1\c x)\1\cc\ ,\qquad
\PR(x):=\1\c\eps(x\1\cc)
\eeq
and introduce the notation $A^L:=\PL(A)$, $A^R:=\PR(A)$ for their
images. The analogue objects in the dual bialgebra $\duA$ will be
denoted by $\duPL,\duPR,\duA^L,$ and $\duA^R$, respectively.
 
Substituting $y=\1$ in Axiom (A.6b) one obtains
immediately the identities
\begin{lrarray}
\eps(x\PL(y))=\eps(xy)\eql&
\eps(\PR(x)y)=\eps(xy)\eqr\label{5}\\
\PL\circ\PL=\PL\eql&\PR\circ\PR=\PR\eqr\label{6}
\end{lrarray}
As a first application of the duality principle
take\footnote{In taking the transpose of a statement with
$\PLR$ use the fact that in a WBA\hfill\break
$\bra\varphi,\PL(x)\ket=\bra\du1\c\o\du1\cc,\1\c\o
x\ket\bra\varphi,\1\cc\ket=\bra\duPL(\varphi),x\ket$ and similarly
$\bra\varphi,\PR(x)\ket=\bra\duPR(\varphi),x\ket$.}
the duals of Eqns(\ref{5}a-b),
\[\1\c\o\PL(\1\cc)\ =\ \1\c\o\1\cc\ =\ \PR(\1\c)\o\1\cc\ .\]
Then these are identities in any WBA. It follows that
\beq
\cop(\1)\ \in\ A^R\ \o\ A^L\ .              \label{9}
\eeq
\begin{lem}
The counit defines a non-degenerate bilinear form
$$x^L\in A^L, y^R\in A^R\ \mapsto\ \eps(y^Rx^L)\in K\ .$$
Hence $A^L\cong A^R$ as $K$-spaces.
\end{lem}
\Proof
\beanon
\eps(y^Rx^L)=0\ \forall y^R\in A^R &\Longrightarrow&
x^L=\eps(\1\c x^L)\1\cc=0\\
\eps(y^Rx^L)=0\ \forall x^L\in A^L &\Longrightarrow&
y^R=\1\c\eps(y^R\1\cc)=0
\eeanon
where we used (\ref{9}).\qed
 
Returning to Eqns(\ref{5}a-b) and substituting them into the
definitions (\ref{Pi}) one obtains
\begin{lrarray}
\PL(x\PL(y))=\PL(xy)\eql&
\PR(\PR(x)y)=\PR(xy)\eqr\label{10}
\end{lrarray}
The duals of (2.5a-b),
\begin{lrarray}
\cop(A^L)\subset A\o A^L\eql&
\cop(A^R)\subset A^R\o A\eqr\label{11}
\end{lrarray}
tell us that $A^L$ and $A^R$ are left, respectively right coideals
in the coalgebra $A$. Using Axiom (A.7b) we can obtain explicit
expressions for these coproducts
\beabc
\cop(x^L)&=\eps(\1\c x^L)\1\cc\o\1\ccc=\eps(\1\cp x^L)\1\c\1\ccp
\o\1\cc=&\1\c x^L\o\1\cc\label{12a}\\
\cop(x^R)&=\1\c\o\1\cc\eps(x^R\1\ccc)=\1\c\o\1\cp\1\cc\eps(x^R\1\ccp)
       =&\1\c\o x^R\1\cc\label{12b}
\eeabc
where $x^L$ and $x^R$ are meant to denote arbitrary elements of
$A^L$, resp. $A^R$.
 
\begin{lem}
For all $x\in A$ we have the identities
\beabc
x\c\o\PL(x\cc)&=&\1\c x\o\1\cc\label{13a}\\
\PR(x\c)\o x\cc&=&\1\c\o x\1\cc\label{13b}
\eeabc
\end{lem}
\Proof Using Axiom (A.7b) one obtains
\beanon
x\c\o\eps(\1\c x\cc)\1\cc&=&\1\cp x\c\eps(\1\c\1\ccp x\cc)\o\1\cc
   =\1\c x\c\eps(\1\cc x\cc)\o\1\ccc=\\
   &=&\1\c x\o\1\cc\\
\1\c\eps(x\c\1\cc)\o x\cc&=&\1\c\o\eps(x\c\1\cp\1\cc)x\cc\1\ccp
   =\1\c\o\eps(x\c\1\cc)x\cc\1\ccc=\\
   &=&\1\c\o x\1\cc
\eeanon
\qed
 
As a consequence we obtain the dual statements
\begin{lrarray}
x\PL(y)=\eps(x\c y)x\cc\eql&
\PR(x)y=y\c\eps(x y\cc)\ .\eqr\label{14}
\end{lrarray}
 
\begin{prop}
Let $A$ be a WBA. Then $A^L$ and $A^R$ are subalgebras of $A$
containing $\1$ and
\beq                       \label{LR=RL}
x^Ly^R\ =\ y^Rx^L\qquad\mbox{for all }x^L\in A^L\ \mbox{and }
           y^R\in A^R\ .
\eeq
\end{prop}
\Proof Eqns(\ref{13a}-b) imply the relations
\beabc
\1\c\1\cp\o\1\cc\o\1\ccp&=&\1\c\o\PL(\1\cc)\o\1\ccc\label{15a}\\
\1\c\o\1\cp\o\1\cc\1\ccp&=&\1\c\o\PR(\1\cc)\o\1\ccc\label{15b}
\eeabc
Now either Axiom (A.7a) or (A.7b) show that on the RHS of
(\ref{15a}) the first tensor factor belongs to $A^R$ and on the
RHS of (\ref{15b}) the last factor belongs to $A^L$. This is
sufficient for $A^R$, respectively $A^L$ to be closed under
multiplication. Hence they are algebras. Obviously $\1\in A^L\cap
A^R$ since $\PL(\1)=\1=\PR(\1)$. In order to see commutativity of
left and right elements just compare Axioms (A.7a) and (A.7b).
\qed
 
As the duals of the statements that $A^L$ and $A^R$ are
subalgebras we obtain that $\Ker\PL$ and $\Ker\PR$ are coideals of
the coalgebra $A$, i.e.
\bea
\cop(\Ker\sqcap^C)&\subset& A\o\Ker\sqcap^C\ +\ \Ker\sqcap^C\o A\
,\\
\eps(\Ker\sqcap^C)&=&0\ ,\qquad\qquad C=L,R\ .\nonumber
\eea
On the other hand, being the annihilator of the left coideal
$\duA^L$, $\Ker\PL$ is a left ideal of the algebra $A$ and
similarly, $\Ker\PR$ is a right ideal.
\begin{lem}
Consider $A^L$ and $A$ as left $A^L$-modules by left
multiplication. Then $\PL\colon$ $A\to A^L$ is a left $A^L$-module
map. Analogously, $\PR\colon A\to A^R$ is a right $A^R$-module
map. That is to say
\beabc
\PL(\PL(x)y)&=&\PL(x)\PL(y)\label{17a}\\
\PR(x\PR(y))&=&\PR(x)\PR(y)\label{17b}
\eeabc
hold true for all $x,y\in A$.
\end{lem}
\Proof At first use the definition of $\PLR$, then
Eqn(\ref{5}a-b), and finally Eqn(\ref{12a}-b):
\beanon
\PL(\PL(x)y)&=\eps(\1\c\PL(x)y)\1\cc=\eps(\1\c\PL(x)\PL(y))\1\cc
             =&\PL(x)\PL(y)\\
\PR(x\PR(y))&=\1\c\eps(x\PR(y)\1\cc)=\1\c\eps(\PR(x)\PR(y)\1\cc)
             =&\PR(x)\PR(y)
\eeanon
\qed
 
Our next assertion about WBA-s establishes a canonical
isomorphism between the left (right) subalgebra of $A$ and the
right (left) subalgebra of $\duA$. Since the existence of a common
non-trivial subalgebra of $A$ and $\duA$ for Hopf algebras is by
far not typical, this result is the first hint towards the
fundamental role $A^L$ and $A^R$ will play in the theory of WHAs.
 
In order to formulate the statement we introduce the Sweedler
arrow notation
\beq                                            \label{Sw arrow}
x\la\varphi:=\varphi\c\bra\varphi\cc,x\ket\ ,\qquad
  \varphi\ra x:=\bra\varphi\c,x\ket\varphi\cc\ .
\eeq
Since $A$ is the dual WBA of $\duA$, the Sweedler arrows
$\varphi\la x$ and $x\ra\varphi$ are also defined.
\begin{lem} \label{kappa}
The map $\kappa_A^L\colon x^L\mapsto (x^L\la\du1)$ is an algebra
isomorphism from $A^L$ onto $\duA^R$.
The map $\kappa_A^R\colon x^R\mapsto(\du1\ra x^R)$ is an algebra
isomorphism from $A^R$ onto $\duA^L$. Furthermore, the restriction
of the canonical pairing to $\duA^L\x A^L$, $\duA^R\x A^R$,
$\duA^L\x A^R$, or to $\duA^R\x A^L$ is non-degenerate.
\end{lem}
\Proof Using Eqns(\ref{15a}-b) and the defining properties
$\bra\varphi\ra x,y\ket=\bra\varphi,xy\ket$,\dots etc. of the
Sweedler arrows one can easily verify that
\bea
(x^L\la\du1)(y^L\la\du1)&=&\du1\c\du1\cp\bra\du1\cc,x^L\ket
                           \bra\du1\ccp,y^L\ket=
              \du1\c\bra\du1\cc,x^L\ket\bra\du1\ccc,y^L\ket=\nn
              &=&x^Ly^L\la\du1\\
(\du1\ra x^R)(\du1\ra y^R)&=&\bra\du1\c,x^R\ket\bra\du1\cp,y^L\ket
                             \du1\cc\du1\ccp=
              \bra\du1\c,x^R\ket\bra\du1\cc,y^R\ket\du1\ccc\nn
              &=&\du1\ra x^Ry^R\\
(\du1\ra x^R)\la \1&=&\1\c\bra\du1\ra x^R,\1\cc\ket=
                     \1\c\eps(x^R\1\cc)=x^R\\
\1\ra(x^L\la\du1)&=&\bra\1\c,x^L\la\du1\ket\1\cc=
                     \eps(\1\c x^L)\1\cc=x^L
\eea
Thus $\kappa_A^L$ ($\kappa_A^R$) is an algebra map with inverse
$\kappa_{\duA}^R$ ($\kappa_{\duA}^L$). As for the non-degeneracy
\beanon
\bra\varphi^R,x^L\ket=0\ \forall\varphi^R&\Longrightarrow&
     x^L=\bra\du1,x^L\1\c)\1\cc=\bra\1\c\la\du1,x^L\ket\1\cc=0\\
\bra\varphi^L,x^L\ket=0\ \forall\varphi^L&\Longrightarrow&
     x^L=\bra\du1,\1\c x^L\ket\1\cc=\bra\du1\ra\1\c,x^L\ket\1\cc=0
\eeanon
and the transpose of these prove the claim. \qed

If $\{b_i\}$ is a $K$-basis of $A$ and $\{\beta^i\}\subset\duA$ is
its dual basis, $\langle\beta^i,b_j\rangle=\delta_{ij}$, then
\beabc
\sum_i\PL(b_i)\o\beta^i&=&\sum_ib_i\o\duPL(\beta^i)\ =\
\1\ra\du1\c\o\du1\cc \label{piba}\\
\sum_i\PR(b_i)\o\beta^i&=&\sum_i b_i\o\duPR(\beta^i)\ =\
\1\c\o\1\cc\la\du1 \label{pibb}
\eeabc
This can be easily seen by pairing both hand sides of any of these
equations with $\varphi\o x$ and apply the definitions (\ref{Pi}).
 
The four arrow identities of the next Scholium will be
frequently used in later computations.
\begin{scho} \label{sch: arrow}
Let $A$ be a WBA. Then for all $\varphi\in\duA$, $x^L\in A^L$, and
$x^R\in A^R$
\beabc \label{44-45}
x^L\la\varphi&=&(x^L\la\du1)\varphi\\
\varphi\ra x^R&=&\varphi(\du1\ra x^R)\ .
\eeabc
\beabc \label{46-47}
\varphi\ra x^L&=&(\du1\ra x^L)\varphi\\
x^R\la\varphi&=&\varphi(x^R\la\du1)\ .
\eeabc
\end{scho}

\subsection{Weak Hopf algebras}
 
In this Subsection we
will show how the existence of an antipode relates $\PL$, $A^L$
with $\PR$, $A^R$ and derive the expected properties
of $S$ that have been axioms in earlier formulations.
The two most important results will be invertibility of the
antipode and separability of the algebras $A^L$ and $A^R$.
Let us start with the question of uniqueness of the antipode.
\begin{lem} \label{lem: uniq S}
The unit, the counit, and the antipode, if exist, are unique. I.e.
if \break
$(A,\mu,u,\cop,\eps,S)$ and $(A,\mu,u',\cop,\eps',S')$ are both
weak Hopf algebras then $u'=u$, $\eps'=\eps$,\break and $S'=S$.
\end{lem}
\Proof The uniqueness of the unit and the counit are obvious.
Therefore $\PL$ and $\PR$ are common in these two WHA-s. In order
to prove $S'=S$ introduce the convolution product
\beq
(f\diamond g)(x)\ :=\ f(x\c)g(x\cc)\ ,\qquad x\in A
\eeq
on functions $f,g\in\Hom_K(A,A)$. This is an associative operation
in terms of which the antipode axioms take the form
\[\id\diamond S=\PL,\qquad S\diamond\id\diamond S=S,\qquad
  S\diamond\id=\PR\ .\]
Now $S'$ satisfies the same equations with the same $\PL,\PR$,
therefore
\[
S'=S'\diamond\id\diamond S'=S'\diamond\PL=S'\diamond\id\diamond S
  =\PR\diamond S=S\diamond\id\diamond S=S\ .\]
\qed
 
As a preparation for the Theorem below notice that the
definitions (\ref{Pi}) have counterparts involving the antipode:
\begin{lrarray}
\PL(x)=\eps(S(x)\1\c)\1\cc\eql &
\PR(x)=\1\c\eps(\1\cc S(x))\eqr\label{18}
\end{lrarray}
As a matter of fact
\beanon
\PL(x)&=&\eps(\1\c\PL(x))\1\cc=\eps(\PL(x)\1\c)\1\cc=
       \eps(x\c S(x\cc)\1\c)\1\cc=\\
      &=&\eps(\PR(x\c)S(x\cc)\1\c)\1\cc=\eps(S(x\c)x\cc S(x\ccc)
        \1\c)\1\cc=\\
      &=&\eps(S(x)\1\c)\1\cc
\eeanon
where in the subsequent equations (\ref{Pi}), (\ref{LR=RL}),
(A.8a), (\ref{5}b), (A.8b), and finally (A.9) have been used.
Eqn(2.23b) can be proven analogously. As the duals of
(2.23a-b) we have automatically the identities
\begin{lrarray}
\PL(x)=S(\1\c)\eps(\1\cc x)\eql&
\PR(x)=\eps(x\1\c)S(\1\cc).\eqr \label{19}
\end{lrarray}
\begin{lem}
In a WHA $A$ the following identities hold:
\beabc
\PL\circ S\ =&\PL\circ\PR&=\ S\circ\PR\label{20a}\\
\PR\circ S\ =&\PR\circ\PL&=\ S\circ\PL\label{20b}
\eeabc
\end{lem}
\Proof It is sufficient to prove the first equalities in
(\ref{20a}) and (\ref{20b}) because the second ones then follow by
duality.
\beanon
\PL\circ S(x)&=&\eps(\1\c S(x))\1\cc=\eps(\1\c S(x\c)x\cc
               S(x\ccc))\1\cc=\\
             &=&\eps(\1\c S(x\c)\PL(x\cc))\1\cc=\eps(\1\c S(x\c)
               x\cc)\1\cc\\
             &=&\PL\circ\PR(x).
\eeanon
In a similar way one can verify $\PR\circ S=\PR\circ\PL$.\qed
 
The above Lemma implies that $S(A^R)\subset A^L$ and
$S(A^L)\subset A^R$. On the other hand Eqns(2.24a-b) say that
$A^L\subset S(A^R)$ and $A^R\subset S(A^L)$. Therefore the
antipode maps $A^L$ onto $A^R$ bijectively and maps $A^R$ onto
$A^L$ bijectively.
 
\begin{thm}   \label{thm: S}
Let $A$ be a WHA. Then the antipode is antimultiplicative and anti
comultiplicative,
\bea
S(xy)&=&S(y)S(x)\qquad\qquad x,y\in A\ ,\label{21}\\
S(x)\c\o S(x)\cc&=&S(x\cc)\o S(x\c)\qquad\qquad x\in A\
                                                   ,\label{22}
\eea
and the restrictions $S|_{A^L}$ and $S|_{A^R}$ are bijections such
that
\beq
S(A^L)\ =\ A^R\ ,\qquad\qquad S(A^R)\ =\ A^L\ .      \label{23}
\eeq
The unit and the counit are $S$-invariant,
\begin{lrarray}
S(\1)\ =\ \1\ ,\eql&
\eps\circ S\ =\ \eps\ .\eqr\label{24}
\end{lrarray}
Furthermore $S\colon A\to A$ is invertible.
\end{thm}
\Proof We have already shown (\ref{23}). (\ref{22}) is the dual of
(\ref{21}) and (2.29a) is the dual of (2.29b).
\beanon
S(xy)&=&S(x\c y\c)x\cc y\cc S(x\ccc y\ccc)=S(x\c y\c)
                   \PL(x\cc\PL(y\cc))=\\
&=&S(x\c y\c)x\cc\PL(y\cc)S(x\ccc)=\PR(\PR(x\c)
                    y\c)S(y\cc)S(x\cc)=\\
&=&S(y\c)\PR(x\c)y\cc S(y\ccc)S(x\ccc)=S(y\c)y\cc S(y\ccc)
                    S(x\c)x\cc S(x\ccc)=\\
&=&S(y)S(x)\ .
\eeanon
Next we prove (2.29b).
\beanon
\eps(S(x))&=&\eps(S(x\c)x\cc S(x\ccc))=\eps(S(x\c)\PL(x\cc))=
            \eps(S(x\c)x\cc)=\eps(\PR(x))=\\
          &=&\eps(x)\ .
\eeanon
In order to prove invertibility of $S$ notice that the descending
chain $A\supset S(A)\supset S^2(A)\supset\dots$ of WHAs all
contain $\1$ by (2.29a). This implies the existence
of $n\in\N$ such that
\[
\1\ \in\ S^{n+1}(A)\ =\ S^n(A)\ \subset\ S^{n-1}(A)\ .
\]
We want to show that this implies $S^n(A)=S^{n-1}(A)$.
Replacing $A$ by $S^{n-1}(A)$ it is
therefore enough to prove invertibility of $S$
under the additional assumption $S^2(A)=S(A)$, implying
\[
\Ker S\ \cap\ S(A)\ =\ 0\ .
\]
In this case let $\bar S:=S|_{S(A)}$, then $\bar S\colon S(A)\to
S(A)$ is bijective and
\[
P_S\ :=\ {\bar S}^{-1}\circ S\ \colon\ A\to S(A)
\]
is a multiplicative idempotent satisfying
\[
P_S(xS(y))\ =\ P_S(x)S(y)\,,\qquad x,y\in A\ .
\]
By (\ref{23}) $A^{L,R}\subset S(A)$.
Now taking into account the identity $x=x_{(1)}
S(x_{(2)})x_{(3)}$,
which follows directly from axioms (A.8a) and (A.4), then
using also $P_S(\1)=\1$ we have
\beanon
P_S(x)&=&P_S( x_{(1)} S(x_{(2)})x_{(3)})=P_S(x_{(1)})
S(x_{(2)})x_{(3)}= P_S(x_{(1)} S(x_{(2)}))x_{(3)}=\\
&=&P_S(\1)x_{(1)} S(x_{(2)})x_{(3)}=x,
\eeanon
so $\Ker P_S=\Ker S=0$.\qed
 
We are now able to derive (versions of) the original antipode
axioms of \cite{BSz,Sz}:
\beabc
x\c\o x\cc S(x\ccc)&=&\1\c x\o\1\cc\\
S(x\c)x\cc\o x\ccc&=&\1\c\o x\1\cc  \label{eq: 2.29b}\\
x\c\o S(x\cc)x\ccc&=&x\1\c\o S(\1\cc)\\
x\c S(x\cc)\o x\ccc&=&S(\1\c)\o\1\cc x
\eeabc
The first two are just rewritings of the bialgebra identities
(\ref{13a}-b). The second two are more delicate.
\beanon
x\c\o S(x\cc)x\ccc&=&x\c\o\eps(x\cc\1\c)S(\1\cc)
                   =x\c\1\cp\o\eps(x\cc\1\ccp\1\c)S(\1\cc)=\\
                  &=&x\c\1\c\eps(x\cc\1\cc)\o S(\1\ccc)
                   =x\1\c\o S(\1\cc)\\
x\c S(x\cc)\o x\ccc&=&S(\1\c)\eps(\1\cc x\c)\o x\cc
                    =S(\1\c)\eps(\1\cc\1\cp x\c)\o \1\ccp x\cc=\\
                   &=&S(\1\c)\o\eps(\1\cc x\c)\1\ccc x\cc
                    =S(\1\c)\o\1\cc x
\eeanon
The following Proposition also holds, if $A$ is just a WBA, see
\cite{WBA}.
 
\begin{prop} \label{pro:sepa}
Let $A$ be a WHA over $K$. Then $A^L$ and $A^R$ are separable
$K$-algebras, in particular, they are semisimple.
\end{prop}
\Proof Recall that an algebra $A$ is separable if and
only if there exists a $q\in A\o A$ such that $(x\o\1)q=q(\1\o x)$
holds for all $x\in A$ and furthermore $\mu(q)=\1$, where $\mu$
denotes the multiplication map of $A$ \cite{Pierce} . Such a $q$ will be
called
a {\em separable idempotent}\footnote{In fact $q$ is an idempotent
only if considered as an element of $A\o A^{op}$.}. So, our proof
will consist of showing that $q^L=S(\1\c)\o\1\cc\in A^L\o
A^L$ and $q^R=\1\c\o S(\1\cc)\in A^R\o A^R$ are separable
idempotents of $A^L$ and $A^R$, respectively. In fact we will
prove the somewhat more general identities
\beabc \label{separ}
x\c y^R\o x\cc&=&x\c\o x\cc S(y^R)\\
x\c\o y^L x\cc&=&S(y^L)x\c\o x\cc
\eeabc
valid for all $x\in A$ and $y^L\in A^L$, $y^R\in A^R$.
Pairing the LHS of (\ref{separ}) with $\varphi\o\psi$, we obtain
\beanon
\bra\varphi\o\psi,LHS\ket&=&\bra\varphi(y^R\la\du1),x\c\ket
\bra\psi,x\cc\ket=\bra\varphi,x\c\ket\bra(S(y^R)\la\du1)\psi,x\cc\ket=\\
&=&\bra\varphi\o\psi,RHS\ket\ .
\eeanon
The proof of (2.31b) is simply the mirror image of the
above argument. \qed

\subsection{The "trivial" representation}
 
Since the counit of a WHA is in general not an algebra map, weak
Hopf algebras may be lacking of any 1-dimensional representation.
Nevertheless the axioms ensure that any WHA $A$ has a
distinguished representation providing a unit object for the
(relaxed) monoidal category of left $A$-modules. We shall discuss
this category in detail in \cite{BNSzII}. Now we concentrate
only on the properties of this representation. We remark that the
trivial representation exists already in WBA's \cite{WBA} and
therefore the use of the antipode in this subsection is not
obligatory.
 
Since the algebras $A^{L/R}$ occur on the right hand side of
Axioms (A.8a-b) where in ordinary Hopf algebras the trivial
representation stands, one expects that the "trivial
representation" of WHAs must be a non-trivial representation
acting on either one of the algebras $A^{L/R}$ or $\duA^{L/R}$.
 
\begin{lem}                     \label{lem: trivials}
The following left $A$-modules are isomorphic.
$$\begin{array}{rccl}
_A\duA^R\ ::\ &\mbox{the vector space }\duA^R&\mbox{with action }&
               x\cdot\varphi^R:=x\la\varphi^R\\
_A\duA^L\ ::\ &\mbox{the vector space }\duA^L&\mbox{with action }&
               x\cdot\varphi^L:=\varphi^L\ra S(x)\\
_AA^L\ ::\ &\mbox{the vector space }A^L&\mbox{with action }&
               x\cdot y^L:=\PL(xy^L)\\
_AA^R\ ::\ &\mbox{the vector space }A^R&\mbox{with action }&
               x\cdot y^R:=\PR(y^RS(x))
\end{array}$$
\end{lem}
\Proof $\duS\colon\duA^L\to \duA^R$ is an isomorphism of vector
spaces and $\duS(\varphi\ra S(x))=x\la\duS(\varphi)$ is a general
WHA identity. This proves the isomorphism of the first two
$A$-modules. Similarly, $S\colon A^L\to A^R$ is an isomorphism
of vector spaces and $S(\PL(xy))=\PR(S(y)S(x))$ is a WHA identity.
This proves the isomorphism of the last two $A$-modules.
 
In order to show the isomorphism of $_A\duA^R$ with $_AA^L$
consider the bijection $B\colon\duA^R\to A^L$,
$B(\varphi^R):=\1\ra\varphi^R$. Then
\beanon
B(x\la\varphi^R)&=&\1\ra(x\la\varphi^R)=\bra\1\c
x,\varphi^R\ket\1\cc=\\
&=&\bra\1\c(x\ra\varphi^R),\du1\ket\1\cc=\PL(x\ra\varphi^R)=
   \PL(x(\1\ra\varphi^R))=\\
&=&\PL(xB(\varphi^R))
\eeanon
hence $B$ is a left $A$-module map. Here, in the last-but-one
equality we have used one of the four arrow identities
of Scholium \ref{sch: arrow}.\qed
\begin{defi}                         \label{def: triv}
By the trivial representation of the WHA $A$ we mean the cyclic
left $A$-module $\V:=\,_A\duA^R$ with $A$-action $\triv\colon A\to
\End_K \duA^R$, $\triv(x)\varphi:=x\la\varphi$.
\end{defi}
The third and fourth $A$-modules of the above Lemma
demonstrate that the restriction of the trivial representation to
$A^L$ ($A^R$) is equivalent to its left regular representation,
hence faithful. This is one of the instances where $A^{L/R}$
appears in the role of a ground "field".
 
Later we will need the following strengthening of Lemma
\ref{kappa}.
\begin{lem}  \label{Z}
Let $A$ be a WHA and introduce the notation $Z^L:=A^L\cap\Center
A$, $Z^R:=A^R\cap\Center A$, and $Z:=A^L\cap A^R$. Then the
isomorphism (of algebras) $\kappa_A^L\colon A^L\to\duA^R$
restricts to an isomorphism $Z^L\to\duZ$ and the isomorphism
$\kappa_A^R\colon A^R\to\duA^L$ restricts to the isomorphism
$Z^R\to\duZ$. Therefore
\beanon
Z^L\la\du1\ =&\duZ&=\ \du1\ra Z^R\\
Z\la\du1=\duZ^R\,,&\ \ &\duZ^L=\du1\ra Z\ .
\eeanon
The two isomorphisms have a common restriction to the {\em
hypercenter} $\Hypercenter A:=Z^L\cap Z^R$ and yields an
isomorphism $\Hypercenter A\to\Hypercenter \duA$.
\end{lem}
\Proof Notice that for $c\in\Center A$ $\du1\ra c=c\la\du1$.
Therefore $x^L\in Z^L\Rightarrow x^L\la\du1=\du1\ra x^L\in \duZ$.
This proves $\kappa_A^L(Z^L)\subset\duZ$.
 
If $z\in Z$ then $(z\la\du1)\varphi=z\la\varphi$
by (\ref{44-45}) and $z\la\varphi=\varphi(z\la\du1)$ by
(2.21b). Hence $z\la\du1$ is central. This proves
$\kappa_A^L(Z)\subset\duZ^R$.
 
Since $(\kappa_A^L)^{-1}=\kappa_{\duA}^R$, the analogue inclusions
$\kappa_A^R(Z^R)\subset \duZ$ and $\kappa_A^R(Z)\subset \duZ^L$
complete the proof.\qed
 
The unusual feature of the trivial representation of WHA-s is that
it can be decomposable. But this can occur only if the left and
right
subalgebras of the dual have non-trivial intersection as the next
Proposition claims.
 
\begin{prop} \label{prop: End trivi}
Let $A$ be a WHA, let $(\V,\triv)$ be its trivial
representation as in Definition \ref{def: triv}. Then
\beq
\End\V\ =\ \triv(Z^L)\ =\ \triv(Z^R)\ ,
\eeq
where $\End\V$ denotes the algebra of $A$-module endomorphisms
of $\V$.
\end{prop}
\Proof Let $T\in\End\V$ then $T(x\la\du1)=x\la T(\du1)$, for $x\in
A$, in particular
\beanon
T(x^L\la\du1)&=& x^L\la T(\du1)\ =\ (x^L\la\du1)T(\du1)\\
T(x^L\la\du1)&=& T(S^{-1}(x^L)\la\du1)=S^{-1}(x^L)\la T(\du1)\ =\
               T(\du1)(x^L\la\du1)
\eeanon
where we have made use of Eqns (\ref{44-45}) and (2.21b).
Since by Lemma \ref{kappa} $A^L\la\du1=\duA^R$,
$\zeta:=T(\du1)\in\Center\duA^R$ and $T(\varphi^R)=\varphi\zeta$.
Thus $x\la\zeta=T(x\la\du1)=(x\la\du1)\zeta$ holds for all $x\in
A$. It follows that
$$
\bra\duPL(\zeta),x\ket=\bra\zeta,\1\ra(x\la\du1)\ket=
\bra(x\la\du1)\zeta,\1\ket=\bra x\la\zeta,\1\ket=\bra\zeta,x\ket\
, $$
i.e. $\zeta\in\duA^L\cap\duA^R\equiv \duZ$. Now by Lemma \ref{Z}
there exists a $z^L\in Z^L$ such that $\zeta=z^L\la\du1$. We can
conclude that
$$
T(\varphi^R)=\zeta\varphi^R=(z^L\la\du1)\varphi^R=z^L\la\varphi^R=
             \triv(z^L)\varphi^R\ ,
$$
i.e. $T=\triv(z^L)$. This proves $\End\V\subset\triv(Z^L)$.
The opposite inclusion is trivial since $\triv(Z^L)\subset
\Center(\triv(A))$. This finishes the proof of
$\End\V=\triv(Z^L)$.
 
Showing the other statement $\End\V=\triv(Z^R)$ one proceeds as
above but chooses a $z^R\in Z^R$ such that $\zeta=\du1\ra z^R$.
Then
\beanon
T(\varphi^R)&=&\varphi^R(\du1\ra z^R)=\varphi^R(z^R\la\du1)=
             z^R\la\varphi^R=\\
            &=&\triv(z^R)\varphi^R
\eeanon
completes the proof.\qed
 
Notice that the above Proposition does not imply that the trivial
$A$-module is semisimple. It does imply, however, that $\V$ has a
decomposition $\V\cong\oplus_{\nu}\,V_\nu$ into indecomposable
$A$-modules in which the indecomposables are disjoint, i.e.
$\Hom(V_\mu,V_\nu)=0$ for all $\mu\neq\nu$.
 
\begin{defi}
If $Z^L=K\1$, or equivalently, if the trivial representation is
indecomposable then the WHA is called {\em pure}.
\end{defi}
The name "pure" comes from the $C^*$-setting when the trivial
representation arises from the positive linear functional $\eps$
by the GNS construction. Thus $A$ is pure iff $\eps$ is pure.
 
Nota bene pureness is not a selfdual notion, duals of pure WHA-s
may not be pure. Clearly, $A$ is pure iff $Z^L\cong Z^R$ is
trivial but $\duA$ is pure iff $Z$ is trivial.

\sec{Weak Hopf Modules and Integral Theory}
 
As in Hopf algebras so in weak Hopf algebras the integrals play
a decisive role in the structure analysis of these algebras.
Using integrals we can formulate conditions for the algebra to be
Frobenius, symmetric, or semisimple, and study questions related
to innerness of $S^2$ or $S^4$. Furthermore we will be able to
characterize those WHAs that have Haar measures. In deriving the
basic properties of integrals the weak generalization of the
the Fundamental Theorem of Hopf modules is very useful. Unfortunately,
it seems to be less powerful than in Hopf algebra theory
(cf. \cite{Montgomery}) where
it implies the existence of non-degenerate integrals. It is an
open problem yet whether all WHAs are Frobenius algebras. We can
prove, however, that all of them are quasi-Frobenius algebras.

\subsection{Integrals in weak Hopf algebras}
The following definition provides the weak Hopf generalization of
the well known notion of integrals in a Hopf algebra \cite{Sweedler}.
 
\begin{defi} \label{def: integral}
A left (right) integral in a weak Hopf algebra $A$ is an
element $l\in A$ {\rm (}$r\in A${\rm )} satisfying
\beq                                                 \label{l (r)}
xl\ =\ \PL(x)l\qquad \left(\ rx\ =\ r\PR(x)\ \right)
\eeq
for all $x\in A$. The space of left (right) integrals in $A$ is
denoted by $\IL$ $(\IR)$. Elements of $\I:=\IL\cap\IR$ are called
two-sided integrals.
A left or right integral in $A$ is called non-degenerate if it
defines a non-degenerate functional on $\duA$. $l\in \IL$ is
called normalized if $\PL(l)=\1$, $r\in \IR$ is called
normalized if $\PR(r)=\1$.
\end{defi}
 
Some equivalent formulations of left (right) integrals are
gathered in the next
\begin{lem}     \label{scho: IL}
Let $A$ be a weak Hopf algebra. Then
the following statements for an element $l\in A$ are equivalent:
\begin{description}
\item[a)] $l\in \IL$
\item[b)] $l\c\o x l\cc=S(x)l\c\o l\cc$ for all $x\in A$
\item[c)] $l\la\duA\subset \duA^L$
\item[d)] $(\varphi\ra x)\la l=S(x)(\varphi\la l)$ for all
$\varphi\in\duA$ and $x\in A$
\item[e)] $(\mbox{\rm Ker}\PL)l=0$
\item[f)] $S(l)\in\IR$
\end{description}
\end{lem}
\Proof {\bf a) $\Rightarrow$ b)}: Using (\ref{eq: 2.29b}) and
(\ref{12a}) we have $l\c\o xl\cc=[S(x\c)\o\1]\cop(x\cc l)=$
$[S(x\c)x\cc S(x\ccc)\o\1]\cop(l)=S(x)l\c\o l\cc$.
{\bf b) $\Rightarrow$ a)}: $xl=x\c l\c\eps(x\cc l\cc)=$\break
$x\c S(x\cc)l\c \eps(l\cc)=\PL(x)l$.
{\bf a) $\Leftrightarrow$ c)}: For an $l\in A$ the equation
$\langle l\la\varphi,x\rangle=\langle\duPL(l\la\varphi),x\rangle$
is clearly equivalent to the equation
$\langle\varphi,xl\rangle=\langle\varphi,\PL(x)l\rangle$.
{\bf b) $\Leftrightarrow$ d)}: By pairing the 2nd tensor factor of
{\bf b)} with an arbitrary $\varphi\in\duA$.
{\bf a) $\Rightarrow$ e)}: is obvious.
{\bf e) $\Rightarrow$ a)}: $xl=[x-\PL(x)]l+\PL(x)l=\PL(x)l$.
{\bf f) $\Leftrightarrow$ a)}: This follows by applying $S$ to
(\ref{l (r)}). \qed

Definition \ref{def: integral} as well as Lemma \ref{scho:
IL} provide rather technical
characterizations of integrals. The next argument sheds some light
on their real nature. Consider the left $A$-module map $\eps_R$
from the left
regular $A$-module to the trivial $A$-module given by acting
with the trivial representation on the cyclic vector $\du1$:
\bea      \label{eq: vartheta}
\eps_R\ \colon&_AA& \to\ _A\duA^R\nn
&x& \mapsto\ (x\la\du1)\ .
\eea
The existence of this (non-zero) map shows that
$\Hom(_AA,\,_A\duA^R)$ is non-zero. However,
there is in general no guarantee that $\Hom(_A\duA^R,\,_AA)$
is non-zero. Left integrals are precisely the objects that label
the possible homomorphisms of the latter type.
\begin{lem}                           \label{lem:embedtrivi}
Left integrals $l$ in $A$ are in one-to-one correspondence with
left $A$-module homomorphisms $f\colon \,_A\duA^R\to \,_AA$. The
correspondence is given by $f\mapsto f(\du1)\in\LI$. What is more
the above map provides an isomorphism
$\LI_A\cong\Hom(\,_A\duA^R,\,_AA)$ of right $A$-modules. In other
words $\LI_A$ is isomorphic to the $A$-dual of the trivial left
$A$-module.
\end{lem}
\Proof If $f\in\Hom(\,_A\duA^R,\,_AA)$ then
$xf(\du1)=f(x\la\du1)=f(\PL(x)\la\du1)=\PL(x)f(\du1)$, hence
$f(\du1)\in\LI$. This is obviously a right $A$-module map. It is
invertible since for $l\in\LI$ the map $f_l\colon\duA^R\to A$,
$f_l(\varphi^R):=(\1\ra\varphi^R)l$ is a left $A$-module map
and satisfies $f_l(\du1)=l$. \qed
 
The identification of $\LI$ with $\Hom(\,_A\duA^R,\,_AA)$ yields
an $A$-valued bilinear form $_A\duA^R\x\LI_A\to A$ given by
evaluation, $(\varphi^R,l)\mapsto f_l(\varphi^R)$.
Replacing $_A\duA^R$ with $_AA^L$ using the isomorphism of Lemma
\ref{lem: trivials} we obtain that this bilinear form is nothing
but multiplication in $A$
\beq                       \label{eq: AL IL}
_AA^L\ \x\ \LI_A\ \to\ _AA_A\ ,\qquad (x^L, l)\ \mapsto x^L \, l
\eeq
and it is an $A$-$A$ bimodule map. We claim that (\ref{eq: AL IL})
is a non-degenerate bilinear form. From one side, $x^Ll=0\
\forall x^L\in A^L\ \Rightarrow\ l=0$, this is trivial. From the
other side we will be able to prove this after having established
that WHAs are quasi-Frobenius algebras in Thm. \ref{th:q-Frob}. As
a matter of fact by Theorem 61.2 of \cite{Curtis-Reiner} the left
annihilator of the right annihilator of the left ideal $\Ker\PL$
is $\Ker\PL$ itself. Now by Lemma \ref{scho: IL}.e) the right
annihilator of $\Ker\PL$ is just $\LI$. Thus $x^Ll=0\ \forall
l\in\LI\ \Rightarrow\ x^L=0$ follows.
 
Now we turn to an other characterization of left integrals
that is related to conditional expectations.
Notice at first that if $\lambda\in \duIL$ then
the map $E_{\lambda}\colon x\mapsto \lambda\la x$ is an
$A^L$-$A^L$-bimodule map from $A$ into $A^L$ commuting with the
right $\hat A$-action on $A$. In fact, all such
maps arise from a left integral, as the following Lemma shows.
\begin{lem}
The left integrals $\lambda\in\duIL$ are in one-to-one
correspondence with right $\duA$-module maps
$E\in\Hom(A_{\duA},A^L_{\duA})$ via
\beanon
\lambda&\mapsto& E_{\lambda}\\
E&\mapsto&\eps\circ E\ .
\eeanon
\end{lem}
\Proof If $\lambda$ is a left integral then $E_{\lambda}$ is a
right $\duA$-module map and maps into $A^L$ by Lemma
\ref{scho: IL}.c).
 
Now let $E\in\Hom(A_{\duA},A^L_{\duA})$. Then
\[
E(x)=\eps(\1\c E(x))\1\cc=\eps\circ E(S^{-1}(\1\c)x)\1\cc
\]
where we used the fact that a right $\duA$-module map is an
$A^L$-$A^L$-bimodule map by (2.20b) and (2.21a).
Hence
\[
E(x)=\eps\circ E(x\ccc)x\c S(x\cc)=\PL(\lambda\la x)
\]
where $\lambda:=\eps\circ E$. It remained to show that $\lambda$
is a left integral.
\[
\bra\varphi\lambda,x\ket=\eps(E(x\ra\varphi))=\eps(E(x)\ra\varphi)
=\bra\duPL(\varphi)\lambda,x\ket
\]
which proves the claim.\qed
 
The characterization of left integrals $\lambda$ as "conditional
expectations" $E_{\lambda}$ provides a link to the theory of
inclusions and "Jones extensions" \cite{NSzW}.
 
The properties of the normalized and the non-degenerate
left integrals will be discussed in later Subsections. Here we
only remark that $\lambda$ is non-degenerate iff $E_{\lambda}$ is
non-degenerate and $\lambda$ is normalized iff $E_{\lambda}$ is
unital.
 
There are two twisting operations $A\mapsto A^{op}$ and $A\mapsto
A_{cop}$ that produce WHAs from WHAs. In the first one the
multiplication $\mu$ is replaced with opposite multiplication
$\mu^{op}(x,y)=\mu(y,x)$ while in the second the coproduct is
replaced by $\Delta^{op}(x)=x\cc\o x\c$. In both cases the
antipode is replaced by $S^{-1}$. The left and right
subalgebras$/$integrals and the dual WHAs of the resulting four
twisted
versions of a WHA $A$ are related to those of $A$ as follows.
 
\[\begin{array}{|l|ccccccc|} \hline
& \PL & \PR & A^L & A^R & \LI & \RI &\duA\\ \hline
A=A(\mu,\cop,S)
&\PL&\PR&A^L&A^R&\LI&\RI&\duA\\
A^{op}=A(\mu^{op},\cop,S^{-1}) &S^{-1}\circ\PR&
S^{-1}\circ\PL&A^L & A^R & \RI & \LI&\duA_{cop}\\
A_{cop}=A(\mu,\cop^{op},S^{-1})
&S^{-1}\circ\PL&S^{-1}\circ\PR&A^R&A^L&\LI&\RI&\duA^{op}\\
A^{op}_{cop}=A(\mu^{op},\cop^{op},S)
&\PR&\PL&A^R&A^L&\RI&\LI&\duA^{op}_{cop}\\ \hline
\end{array}\]
As an application of the table we give here the twisted versions
of the identity of Lemma \ref{scho: IL}.d):
\beabc
(\varphi\ra x)\la l& =& S(x) (\varphi\la l)\\
(x\la\varphi)\la r&=&(\varphi\la r) S^{-1}(x)\\
l\ra(\varphi\ra x)&=& S^{-1}(x)(l\ra\varphi)\\
r\ra(x\la\varphi)&=&(r\ra\varphi) S(x)
\eeabc
for all $x\in A$, $\varphi\in \duA$, $l\in \LI$, and $r\in \RI$.
 
\subsection{Weak Hopf modules}
Let $A$ be a WHA. Recall that a {\em left $A$-module} is a
$K$-linear
space $M$ carrying a left action of the algebra $A$, denoted by
$x\in A,m\in M\mapsto x\cdot m$. A {\em right $A$-module} is a
left module $M$ of the opposite algebra $A^{op}$ with action
denoted by $x\in A,m\in M\mapsto m\cdot x$. Since $A$ is unital,
all modules are assumed to be non-degenerate, i.e. $\1$ acts as
the identity. The left $A$-module $M$ is called {\em faithful} if
$x\cdot m=0$, $\forall m\in M$ implies $x=0$.
 
The $A$-modules know nothing about the coalgebra structure of $A$.
The {\em left $A$-comodules} $M$ in turn are the comodules of the
coalgebra $A$ and carry no information about the algebra structure
of $A$. The left coaction is denoted by $m\mapsto m_{-1}\o m_0\in
A\o M$. One defines the {\em right $A$-comodules} analogously and
denotes the coaction as $m\mapsto m_0\o m_1\in M\o A$.
 
Because of the finite dimensionality of $A$ there is a one-to-one
correspondence between left (right) $A$-coactions on $M$ and right
(left) $\duA$-actions on $M$ given by
\bea
m\cdot\varphi&=\ \bra\varphi,m_{-1}\ket m_0\ ,\qquad
m_{-1}\o m_0\ =&\sum_i b_i\o m\cdot\beta^i\\
\varphi\cdot m&=\ m_0\bra\varphi,m_1\ket\ ,\qquad
m_0\o m_1\ =&\sum_i \beta^i\cdot m\o b_i\ .
\eea
Here $\{b_i\}$ denotes an arbitrary basis of $A$ and $\{\beta^i\}$
is its dual basis: $\bra\beta^i,b_j\ket=\delta_{ij}$.
There are 8 basic examples of $A$ modules with the target space
$M$ being either $A$ itself or its dual $\duA$. These are the
following.
$$\begin{array}{rlcrl}
_AA::& x\cdot y=xy\,,&\qquad&A_A::&y\cdot x=yx\\
^AA::& x\cdot y=yS(x)\,,&\qquad&A^A::&y\cdot x=S(x)y\\
_A\duA::&x\cdot\varphi=x\la\varphi\,,&\qquad&\duA_A::&\varphi\cdot
                                                x=\varphi\ra x\\
^A\duA::&x\cdot\varphi=\varphi\ra S(x)\,,&\qquad&\duA^A::&
                       \varphi\cdot x=S(x)\la\varphi
\end{array}$$
where the Sweedler arrow notation (\ref{Sw arrow}) has been used.
They all are faithful and non-degenerate due to the existence of a
unit and a counit.
To each of the $A$-modules in the above list there is a
corresponding $\duA$-comodule denoted by the same symbol.
This identification is justified also by the fact that $N\subset
M$ is an $A$-submodule if and only if it is an $\duA$-subcomodule.
 
By analogy with our definition of left integrals,
the space of {\em invariants}  of a left
$A$-module $M$ is defined to be the subspace
\beq
\Inv M\ :=\ \{\,m\in M\,|\, x\cdot m=\PL(x)\cdot m,\ \forall x\in
A\,\}
\eeq
By the same methods as in Lemma \ref{lem:embedtrivi},
$\Inv M$ is linearly isomorphic to $\Hom({_A\hat A^R},{_AM})$ via
\beq
\Inv M=\{f(\hat\1)\,|\, f\in\Hom({_A\hat A^R},{_AM})\}.
\eeq
By duality, we define the  {\em coinvariants} of a right
$A$-comodule $M$ as
\beq
\Coinv M\ :=\ \{\,m\in M\,|\, m_0\o m_1=m_0\o\PL(m_1)\}
\eeq
Thus, $m\in\Coinv M \Lra m_0\o m_1\in M\o A^L$ and for a
left $A$-module $M$, the invariants $\Inv M\subset M$  coincide with
$\Coinv M\subset M$ considered as an $\hat A$-comodule.
Similarly, for a right $A$-module (left $A$-comodule) $M$
the invariants (coinvariants) are
\bea
\Inv M &=& \{\,m\in M\,|\,m\cdot x=m\cdot\PR(x),\ \forall x\in A\,\}\\
\Coinv M &=&\{\,m\in M\,|\,m_{-1}\o m_0=\PR(m_{-1})\o m_0\}
.
\eea
Notice that the (co)invariants do not form a
sub(co)module, even not an
$A^{L/R}$-sub\-mod\-ule.
 
\begin{scho} \label{scho3.5}
The invariants of the left (right) regular $A$-module
are precisely the left (right) integrals of $A$:
\[
\Inv _AA\ =\ \IL\ ,\qquad \Inv A_A\ =\ \IR\ .
\]
The invariants of $_A\duA$ and $\duA_A$, on the other hand, yield
the left and right subalgebras, respectively:
\[
\Inv _A\duA\ =\ \duA^L\ ,\qquad \Inv \duA_A\ =\ \duA^R\ .
\]
\end{scho}
Investigating the structure of the mixed modules $_{\duA}\duA^A$
and
$^A\duA_{\duA}$, that incorporate the whole bialgebra structure of
$A$, one arrives to a weak generalization of the notion of Hopf
modules \cite{Abe,Sweedler}.
 
\begin{defi}                      \label{def: WHM}
A {\em right weak Hopf module} (right WHM) over $A$ is a right $A$
module $M$ which is also a right $A$-comodule such that
the compatibility relation
\beq                                   \label{WHM comp}
(m\cdot x)_0\o (m\cdot x)_1\ =\ m_0\cdot x\c\o m_1 x\cc
\eeq
holds for $x\in A,\ m\in M$.
\end{defi}
 
\begin{lem}      \label{lem3.7}
Let $M$ be a right WHM over $A$. Then for all $m\in M$
 
i) $m_0\cdot\PR(m_1) =m$
 
ii) $\Coinv M =\{m\in M\,|\,m_0\o m_1=m\cdot\1\c\o\1\cc\}$
and $\Coinv M$  is a right $A^L$ submodule.
 
iii) $E(m):=m_0\cdot S(m_1)$ provides a projection
$E:M\to\Coinv M$.
\end{lem}
 
{\em Proof:}
i)  Let $M$ be a right WHM over $A$. Since
\beq\label{*}
m_0 \eps(m_1 x)= m_0\cdot \1_{(1)}\eps(m_1\1_{(2)} x)=
m_0\cdot \1_{(1)}\eps(m_1\1_{(2)})\eps(\1_{(3)} x)=
m\cdot \1_{(1)}\eps(\1_{(2)} x)
\eeq
for all $x\in A$, we have $m_0\o {\hat \1}\ra m_1=
m\cdot \1_{(1)} \o {\hat \1}\ra \1_{(2)}$, so that
\beq m_0\cdot \PR(m_1)=
m_0\cdot [({\hat \1}\ra m_1)\la \1]=
m\cdot \1\c [(\du1\ra \1\cc)\la \1]=
m\cdot \1_{(1)}\PR (\1_{(2)})=m.
\eeq
ii) The inclusion $\supset$ follows from (\ref{9}).
Conversely, if $m\in\Coinv M$ then $m_0\o m_1\in M\o A^L$,
implying by (2.3a), (\ref{LR=RL}) and (\ref{*})
$$
m_0\o m_1=m_0\o\eps(m_1\1\c)\1\cc =
m\cdot\1\c\o\eps(\1\cc\1_{(1')})\1_{(2')}=m\cdot\1\c\o\1\cc\,.
$$
Also $\cop({A^L})\subset A\o A^L$ and therefore $\Coinv M$ is a
right $A^L$-submodule.
 
iii) To check that $n:=m_0\cdot S(m_1)$ is a
coinvariant for all $m\in M$ we compute
\beanon
n_0\o n_1&=&(m_0\o m_1)(\cdot (S(m_2)\c\o S(m_2)\cc)=\\
   &=&m_0\cdot S(m_3)\o m_1 S(m_2)=m_0\cdot S(\1\cc m_1)\o
           S(\1\c)=\\
   &=&n\cdot \1\c\o\1\cc .
\eeanon
Since for $n\in \Coinv M$ we have
$n_0\cdot S(n_1)=n\cdot\1\c S(\1\cc)=n$, $E$ is a projection
onto $\Coinv M$.
\qed
 
\begin{exam}\label{ex:rWHM} The right weak Hopf module
$_{\duA}\duA^A$
 
\em As a linear space the module is the dual WHA $\duA$. The right
action and coaction are
\beq
\varphi\cdot x:= S(x)\la \varphi\,,\qquad
\varphi_0\o\varphi_1:=\sum _i\ \beta^i\varphi\o b_i\ .
\eeq
Clearly, the right $A$-coaction is dual to the left
$\hat A$-multiplication and therefore counital and
right-coassociative.
The compatibility condition (\ref{WHM comp}) can be seen as
follows.
\beanon
(\varphi\cdot x)_0\o (\varphi\cdot x)_1
     &=&\sum_i\beta^i(S(x)\la \varphi)\o b_i=\\
     &=&\sum_i\beta^i\left[S(x\c)\la x\cc
               S(x\ccc)\la\varphi\right]\o b_i=\\
     &=&\sum_i\beta^i\left[S(x\c)\la(x\cc\la\du1)\varphi\right]
               \o b_i=\\
     &=&\sum_i\beta^i\left[S(x\cc)x\ccc\la\du1\right]
              \left[S(x\c)\la\varphi\right]\o b_i=\\
     &=&\sum_i\left[S(x\cc)x\ccc\la\beta^i\right]
              \left[S(x\c)\la\varphi\right]\o b_i=\\
&=&\sum_iS(x\c)\la\left[(x\cc\la\beta^i)\varphi\right]
              \o b_i=\\
     &=&\sum_i S(x\c)\la(\beta^i\varphi)\o b_i x\cc=\\
     &=&\varphi_0\cdot x\c\o\varphi_1x\cc\ .
\eeanon
The $A$-coinvariants of this WHM coincide with the $\hat A$-invariants of
the dual left regular $\hat A$-module $_{\hat A}\hat A$ and therefore with
the space of left integrals in $\hat A$ by Scholium \ref{scho3.5},
\beq
\Coinv(_{\duA}\duA^A) = \duIL\,.
\eeq
\end{exam}
 
The Fundamental Theorem of Hopf modules generalizes to the weak
case as follows.
\begin{thm}
Let $A$ be a WHA, $M$ be a right WHM over $A$ and
let $N=\Coinv M$ denote
the set of coinvariants of $M$. Since $N$ is a
right $A^L$ submodule, one can form the $A^L$-module tensor
product $N\amalgo{A^L}A$ and make it into a right WHM by the
definitions
\beabc
(n\o a)\cdot x&:=&n\o ax\\
(n\o a)_0\o (n\o a)_1&:=&(n\o a\c)\o a\cc
\eeabc
where $a,x\in A, n\in N$. Then the map
\beq
\alpha\ \colon\ N\amalgo{A^L}A\ \to\ M\ ,\quad
n\o x\ \mapsto\ n\cdot x
\eeq
is an isomorphism of right WHM-s.
\end{thm}
Recall that an {\em isomorphism of WHM-s} is just a module
isomorphism which is a comodule isomorphism at the same time.
 
\Proof That $\alpha$ is a module map and comodule map is easy to
verify. To construct the inverse define
\bea
\beta&\colon&M\ \to\ N\amalgo{A^L}A\nn
\beta(m)&=&m_0\cdot S(m_1)\o m_2\ \equiv\ E(m_0)\o m_1\ .
\eea
Then $\beta$ is obviously a comodule map.
We show that it is also a module map.
\beanon
\beta(m\cdot x)&=&m_0\cdot x\c S(m_1 x\cc)\o m_2 x\ccc=\\
   &=&m_0\cdot \PL(x\c)S(m_1)\ \o\ m_2 x\cc=\\
   &=&m_0\cdot S(m_1\1\c)\ \o\ m_2\1\cc x=\\
   &=&\beta(m)\cdot x
\eeanon
We are left with showing that on the one hand
\beanon
\beta\circ\alpha(n\o x)&=&\beta(n\cdot x)=\beta(n)\cdot x=
       n\cdot\1\c S(\1\cc)\o \1\ccc x=\\
      &=&n\o\1\c S(\1\cc)\1\ccc x=\\
      &=&n\o x
\eeanon
and on the other hand
\[\alpha\circ\beta(m)\ =\ m_0\cdot S(m_1)m_2\ =\ m\]
where in the last equality Lemma \ref{lem3.7}i) has been used. \qed
 
\bigskip\noindent
Applying this Theorem to the WHM of Example \ref{ex:rWHM}
we obtain the right WHM isomorphism
\beq \label{eq: fundamenta}
_{\duA}\duA^A\ \cong\ \I^L(\duA)\ \amalgo{A^L}\ _{\duA}A_A\ .
\eeq
\begin{coro}                        \label{coro: WHM}
In any WHA $A$ the space of left integrals
$\IL=S(\IR)$ is non-zero and $\duIL$ is the dual of $\IR$
with respect to the restriction of the canonical pairing.
Moreover, choosing a basis $\{\lambda^a\}$ in $\duIL$ and taking
its dual basis $\{r_a\}$ in $\IR$, we have
\bea
\du1&=&\sum_a S(r_a)\la\lambda_a  \label{eq:r-lambda}\\
\1&=&\sum_a r_a\ra\duS(\lambda_a)\  \label{eq:r-lambda-d}
\eea
\end{coro}
\Proof $\duIL\neq 0$ follows from (\ref{eq: fundamenta}).
By inspecting the form of the projection $E: M\to N$
in Example \ref{ex:rWHM} we get a projection
$L\colon A\to\IL$ onto the left integrals,
\beq                                    \label{eq: L}
L(x)\ =\ \sum_i\ \duS^2(\beta^i)\la(b_ix)\ .
\eeq
Therefore the projection to the right integrals is
\beq
R(x)=S\circ L\circ S^{-1}(x)\ =\ \sum_i\ (xb_i)\ra\duS^2(\beta^i)\
.
\eeq
Similar expressions define the projections $\hat L$ and $\hat R$
to the dual integrals. Now it is easy to check that
\beq
\bra \hat L(\varphi), x\ket\ =\ \bra\varphi, R(x)\ket
\eeq
proving the non-degeneracy of the restriction of the canonical
pairing to $\duIL\x\IR$.
 
The dual bases satisfy $\bra\lambda^a,r_b\ket=\delta_{ab}$
therefore
\bea
  \bra \lambda^a\,,\,x S(r_a)\ket &=&
      \bra {\hat L}(\beta^i)\,,\, x S(b_i) \ket =
      \bra S^2(b_j)\la \beta^j \beta^i\,,\, x S(b_i) \ket=
      \bra \beta^k , x S\circ \PR(b_k) \ket=\\
     &=& \bra \1\cc\la \du1 , x S(\1\c) \ket
\eea
where in the last step we used (\ref{pibb}). This proves
(\ref{eq:r-lambda}).
(\ref{eq:r-lambda-d}) is the twisted version in $A^{op}_{cop}$.
\qed

\subsection{Restrictions on the algebraic structure}
The existence of a weak Hopf structure on the $K$-algebra $A$
involves certain restrictions on the algebra $A$, just like in
case of Hopf algebras. In this subsection we show that any WHA
$A$ is quasi-Frobenius, i.e. self-injective. The notions of
semisimple and separable algebras coincide within the class of
WHAs. Moreover, we prove an analogue of Maschke's Theorem which
claims that $A$ is semisimple if and only if it has normalized
left integrals.
 
\begin{thm}     \label{th:q-Frob}
Every weak Hopf algebra over a field $K$ is a
quasi-Frobenius algebra.
\end{thm}
\Proof By Theorem 61.2 of \cite{Curtis-Reiner} it is sufficient
to prove that the left regular $A$-module $_AA$ is injective. By
the Nagao--Nakayama theorem injectivity of a left $A$-module is
equivalent to that it is a direct sum of $K$-duals of principal
indecomposable right $A$-modules.
Since $\duA^A$ is the $K$-dual of $_AA$, we need to
show that $\duA^A$ is a direct sum of principal
indecomposable right $A$-modules, i.e. that $\duA^A$ is
projective.
This in turn is a consequence of the Fundamental Theorem of WHM-s.
 
As a matter of fact we have the right $A$-module isomorphisms
\beq        \label{eq:q-Frob}
\duA^A\ \cong\ \duIL\amalgo{A^L}A_A\ \cong P(\duIL\amalgo{K}A_A)
\eeq
the first of which is the consequence of the Fundamental Theorem
of the right WHM $_{\duA}\duA^A$, the second of which is a rather
simple property of the amalgamated tensor product with
respect to the separable algebra $A^L$. In order to explain the
projection $P$ here we make a digression.
\begin{lem}
Define the map $P\colon \duIL\amalgo{K}A\to\duIL\amalgo{K}A$ by
\beq
P(\lambda\o x)\ :=\ S^2(\1\c)\la\lambda\o\1\cc x\ .
\eeq
Then $P\circ P=P$ and $\Ker P$ coincides with $\Ker\pi$ of the
canonical projection $\pi$ from the free right $A$ module
$\duIL\amalgo{K}A_A$ onto $\duIL\amalgo{A^L}A_A$. Therefore
\beq
P(\duIL\amalgo{K}A_A)\stackrel{\pi|_{{\scriptscriptstyle Im}
P}}{\longrightarrow} \duIL\amalgo{A^L}A_A
\eeq
is an isomorphism of right $A$-modules.
\end{lem}
\Proof The kernel of the canonical projection is
\beq
\Ker\pi\ =\ \mbox{Span}_K\{(\lambda\o x^Ly)-(S(x^L)\la\lambda\o
x)\ |\ \lambda\o y\in\duIL\o A,\ x^L\in A^L\ \}
\eeq
If $\sum_i(\lambda_i\o x_i)\in\Ker\pi$ then obviously $\sum_i
S^2(\1\c)\la\lambda_i\o\1\cc x_i=0$, therefore $\Ker\pi\subset\Ker
P$. Now assume $\sum\lambda_i\o x_i\in\Ker P$. Then
\beanon
\sum_i\lambda_i\o x_i&=&\sum_i\lambda_i\o S(\1\c)\1\cc x_i=\\
&&\sum_i[\lambda_i\o S(\1\c)(\1\cc
x_i)-S^2(\1\c)\la\lambda_i\o\1\cc x_i]\ \in\ \Ker\pi
\eeanon
This proves $\Ker\pi=\Ker P$. That $P$ is a projection and a right
$A$-module map is trivial to verify. Therefore $\pi|_{{\rm
\scriptstyle Im} P}$ is an $A$-module isomorphism.\qed
 
{\it Back to the Proof of Theorem \ref{th:q-Frob}}\,: In virtue of
the above Lemma the amalgamated tensor product
$\duIL\amalgo{A^L}A_A$ is the direct summand of a free $A$-module,
hence projective. By Equation \ref{eq:q-Frob} this is isomorphic
to $\duA^A$. This proves projectivity of $\duA^A$, hence
injectivity of $_AA$.\qed
 
\bigskip\noindent
The equivalence of (c) and (d) of the next Theorem
provides a weak Hopf version of Maschke's Theorem known for Hopf
algebras as well \cite{LS}. Below we denote
$\eps_R(x):=x\la\du1$.
 
\begin{thm}       \label{thm: semi}
The following conditions on a WHA $A$ over $K$ are equivalent:
\begin{description}
\item[a)] $A$ is semisimple.
\item[b)] In the category of left $A$-modules the following exact
sequence is split
$$
0\longrightarrow\Ker\eps_R\longrightarrow\,
_AA\stackrel{\eps_R}{\longrightarrow}\, _A\duA^R\longrightarrow 0
$$
\item[c)] There exists a normalized left integral $l\in A$.
\item[d)] $A$ is a separable $K$-algebra.
\end{description}
\end{thm}
\Proof
(a$\Rightarrow$c):
If $A$ is semisimple, then $\Ker\PL\equiv\Ker\eps_R$ being a
a left ideal there exists $p=p^2\in A$ such that $\Ker\PL=Ap$,
whence $l=\1-p$ is a normalized left integral by Lemma
\ref{scho: IL}.e.
 
(b$\Leftrightarrow$c): Let $F\in\Hom(_A\duA^R,\,_AA)$ be such that
$\eps_R\circ F=\id$. Then $xF(\du1)=F(x\la\du1)=\PL(x)F(\du1)$,
for $x\in A$, therefore $F(\du1)\in\IL$. Moreover,
$\du1=\eps_R(F(\du1))=F(\du1)\la\du1$ implying
$\PL(F(\du1))=\1$.
Conversely, if $l\in\LI$ is a normalized left integral then
$F\in\Hom(_A\duA^R,{_AA})$ given by $F(\eps_R(x)):=xl$ satisfies
$\eps_R\circ F=\id$.
 
(c$\Rightarrow$d): Let $l$ be a normalized left integral. Then
$q=l\c\o S(l\cc)$ is a separating idempotent for $A$. As a matter
of fact $\mu(q)=\1$ follows from the normalization $\PL(l)=\1$
while $(x\o\1)q=q(\1\o x)$ is precisely the left integral property
of Lemma \ref{scho: IL}.b.

(d$\Rightarrow$a):
This is a standard result \cite{Pierce}.
 \qed
 
\subsection{Non-degenerate integrals}
 
Until now we have not been able to decide whether the WHM Theorem
of Subsection 3.2 implies the existence of non-degenerate
integrals, as it does in case of Hopf algebras. In the present
subsection we will show that the existence of non-degenerate
integrals in the WHA $\duA$ is equivalent to the existence of
non-degenerate functionals on $A$, i.e. that $A$ is a Frobenius
algebra. As a byproduct we obtain that the class of Frobenius WHAs
is selfdual.
 
The space $\RI$ of right integrals can be viewed as a $K$-module,
as a left $A^L$-module $_{A^L}\RI$ by left multiplication, and as
a left $A$-module $_A\RI$ since it is a left ideal of $A$. From
the latter point of view $_A\RI$ is the dual of the trivial right
$A$-module, $_A\RI\cong\Hom(\duA^L_A,A_A)$, by a twisted version
of Lemma \ref{lem:embedtrivi}. As a $K$-module $\RI$ has $\duLI$
as its $K$-dual,
$\duLI\cong\Hom(_K\RI,\,_KK)$ the isomorphism being given by the
restriction of the canonical pairing (see Corollary \ref{coro:
WHM}). The
next Lemma shows that $\duLI$ is also the $A^L$-dual of $\RI$ with
right $A^L$-module structure precisely the one needed in Example
\ref{ex:rWHM}, i.e. $\lambda\cdot x^L=S(x^L)\la\lambda$.
\begin{lem}
The $A^L$-valued bilinear form
\beq
(\ ,\ )_{A^L}\ \colon\ \RI\x\duLI\ \to\ A^L\ ,\qquad
(r,\lambda)_{A^L}=\lambda\la r
\eeq
provides an isomorphism of right $A^L$-modules
\beq
(\duLI)^{A^L}\ \to\ \Hom(_{A^L}\RI,\,_{A^L}A^L)\ ,\qquad
\lambda\mapsto\ (r\mapsto(r,\lambda)_{A^L})\ ,
\eeq
i.e. $(\duLI)^{A^L}$ is the $A^L$-dual of $_{A^L}\RI$.
\end{lem}
\Proof At first verify the following properties of the
$A^L$-valued bilinear form.
\beabc
(x^L\cdot r,\lambda)_{A^L}&=&x^L(r,\lambda)_{A^L}\\
(r,\lambda\cdot x^L)_{A^L}&=&(r,\lambda)_{A^L}\,x^L\\
(r,\lambda)_{A^L}=0\ \forall r\in\RI&\Rightarrow&\lambda=0
\eeabc
The first two are simple WHA identities. The third one follows
from the relation $\eps((r,\lambda)_{A^L})=\bra\lambda,r\ket$
and from non-degeneracy of the canonical pairing $\bra\ ,\ \ket$
on $\duLI\x\RI$ (Corollary \ref{coro: WHM}). Now properties (a)
and (b) tell
us that $\lambda\mapsto (.,\lambda)_{A^L}$ is indeed the required
$A^L$-module map and (c) ensures that it is injective. In order to
show that it is surjective it is sufficient to find finite sets of
elements $\{r_a\}$ in $\RI$ and $\{\lambda^a\}$ in $\duLI$ such
that
\beq             \label{fgp}
\sum_a (r,\lambda^a)_{A^L}\ r_a\ =\ r\qquad \forall r\in\RI\ .
\eeq
For if such elements exist then any
$f\in\Hom(_{A^L}\RI,\,_{A^L}A^L)$ can be written as $f=\sum_a
\lambda^a\cdot f(r_a)$. As a matter of fact
\[
f(r)=f(\sum_a (r,\lambda^a)_{A^L}\cdot r_a)=\sum_a
(r,\lambda^a)_{A^L}\ f(r_a)=(r,\sum_a\lambda^a\cdot f(r_a))_{A^L}
\]
for all $r\in\RI$. Now we claim that a pair of dual bases $\{r_a\}$ of
$\RI$ and $\{\lambda^a\}$ of $\duLI$, in the sense of $K$-duality,
i.e. $\bra\lambda^a,r_b\ket=\delta_{ab}$, also satisfies
(\ref{fgp}). As a matter of fact for $\lambda\in\duLI$ we have
\beanon
&&\bra\lambda,\sum_a (\lambda^a\la r)r_a\ket=\eps(\sum_a
(\lambda^a\la r)(\lambda\la r_a))=\sum_a\bra\lambda^a,r(\lambda\la
r_a)\ket=\\
&=&\sum_a\bra\lambda^a, r(S(r_a)\ra \duS^{-1}(\lambda))\ket=
\sum_a\bra\lambda^a\c,r\ket\,\bra\duS^{-1}(\lambda)\lambda^a\cc,
S(r_a)\ket=\\
&=&\sum_a\bra\lambda\lambda^a\c,r\ket\,\bra\lambda^a\cc,S(r_a)\ket=
\bra\lambda\left[\sum_a S(r_a)\la\lambda^a\right],\,r\ket
=\bra\lambda,r\ket
\eeanon
where in the last equality (\ref{eq:r-lambda}) has been used.\qed
 
Notice that Eqn(\ref{fgp}) means that $_{A^L}\RI$ is finitely
generated projective\footnote{Although this is clear from the
fact that $A^L$ is semisimple, constructing the concrete bases
$\lambda^a,\ r_a$ was not in vain since it will help to compute
the commutant in (\ref{eq: EndIR}).}.
Therefore by a general result (see e.g. \cite{Curtis-Reiner new})
\beq
\duLI\amalgo{A^L}\RI\ \cong\ \End\,_{A^L}\RI\ .
\eeq
On the other hand the isomorphism $\alpha$ of the WHM Theorem, if
restricted to $\duLI\amalgo{A^L}\RI$, yields an isomorphism onto
$\duA^L$. Thus we have the composition
\beq
{\cal E}\colon\ \
A^L\longrightarrow\duA^L\stackrel{\alpha^{-1}}{\longrightarrow}
\duLI\amalgo{A^L}\RI\longrightarrow\End\,_{A^L}\RI
\eeq
of isomorphisms. Evaluating it explicitely we obtain
\bea            \label{eq: EndIR}
r\cdot{\cal
E}(x^L)&=&r\cdot\left(\sum_{ij}S^2(b_i)\la(\beta^i\beta^j(\du1\ra
x^L))\ \o\ b_j\right)=r\cdot\left(\sum_j\hat L(\beta^j(\du1\ra
x^L))\o b_j\right)=\nn
&=&r\cdot(\sum_a\lambda^a\o S(x^L)r_a)=\sum_a (r,\lambda^a)_{A^L}
S(x^L)r_a=\nn
&=&S(x^L) r\quad \forall r\in\RI,\ x^L\in A^L\ .
\eea
This proves our next
\begin{prop}
The left modules $_{A^L}\RI$ and $_{A^R}\RI$ are faithful and the
endomorphism algebra of $_{A^L}\RI$ consists of left
multiplications with elements of $A^R$. Therefore
\beq
\End\,_{A^L}\RI\ \cong\ A^L\ ,\qquad\mbox{ as algebras.}
\eeq
\end{prop}
 
The set $\Sec A^L$ of equivalence classes of simple left
$A^L$-modules will be called the {\em sectors} of $A^L$. For $a\in
\Sec A^L$ let $V_a$ be a simple module from the class $a$ and let
$\D_a=\End V_a$ be the corresponding division algebra. Then by the
Wedderburn structure theorem $A^L\cong\oplus_a M_{n_a}(\D_a)$. Let
$m_a$ denote the multiplicity of $V_a$ in the semisimple module
$_{A^L}\RI$. Then $\End\,_{A^L}\RI\cong\oplus_a M_{m_a}(\D_a)$
which is, by the Proposition, isomorphic to $A^L$. This is
possible only if there is a permutation
\beq
\tilde{\empty}\ \colon\ \Sec A^L\to\Sec A^L\,,\quad \mbox{such
that}\ n_{\tilde a}=m_a\ \mbox{and}\ \D_{\tilde a}=\D_a\ .
\eeq
This means that $\RI$, as an $A^L$-$\End\,_{A^L}\RI$ bimodule, can
be identified with a direct sum of matrices,
\beq
_{A^L}\RI\ \cong\ \oplus_a\ \Mat(n_a\x m_a,\, \D_a)\ .
\eeq
This allows us to compute its $K$-dimension and apply the
Cauchy-Schwarz inequality to obtain the bound
\beq
\dim_K\RI\ =\ \sum_a(\dim_K \D_a)\,n_am_a\ \leq\ \sum_a(\dim_K
\D_a)\, n_a^2 \ =\ \dim_K A^L\ .
\eeq
Equality holds here iff $m_a=n_a$, $a\in\Sec A^L$, i.e. iff
$_{A^L}\RI\cong \,_{A^L}A^L$. Now we are ready to prove the
following
\begin{thm}
Let $A$ be a WHA over the field $K$. Then the following conditions
are equivalent.
\begin{description}
\item[i)] $A$ is a Frobenius algebra;
\item[ii)] $\dim_K\,\RI\ =\ \dim_K\,A^L$;
\item[iii)] Non-degenerate integrals exist in $A$;
\item[iv)] $\duA$ is a Frobenius algebra.
\end{description}
\end{thm}
\Proof {\bf (i)$\Rightarrow$(ii)} If $_AA\cong\,_A\duA$ then their
invariants $\I(_AA)=\LI$ and $\I(_A\duA)=\duA^L$, respectively
(see Scholium 3.5), are isomorphic as $K$-spaces.
{\bf (ii)$\Rightarrow$(iii)} As we have seen above the $K$-space
isomorphism of $\RI$ and $A^L$ implies that $_{A^L}\RI$ is
isomorphic to the left regular module $_{A^L}A^L$. Since the
latter is cyclic, there exists a cyclic vector $r\in\,_{A^L}\RI$.
Thus $l:=S(r)$ is cyclic in $(\LI)^{\duA^L}$.  As a matter of
fact $\LI=S(\RI)=S(A^Lr)=lA^R=\duS(\duA^L)\la l$. Now
interchanging the roles of $A$ and $\duA$ in the WHM Theorem
\[
A=\alpha(\LI\amalgo{\duA^L}\duA)=\duS(\duA)\la\LI=\duS(\duA)\la
(\duS(\duA^L)\la l)=\duA\la l\ ,
\]
hence $l$ is a non-degenerate left integral in $A$. {\bf
(iii)$\Rightarrow$(iv)} is obvious since $l$ is a non-degenerate
functional on $\duA$. {\bf (iv)$\Rightarrow$(i)} Repeat the
arguments above from {\bf (i)} to {\bf (iv)} with $A$
replaced by $\duA$. \qed
 
Weak Hopf algebras satisfying any one of the conditions of the
above Theorem will be called {\em Frobenius WHAs}.
Note that since semisimple algebras are Frobenius, in a semisimple WHA there
exist both normalized and non-degenerate integrals, although there
may be no integral sharing both properties.
\footnote{As an example consider $M_2(\Z_2)$, the semisimple algebra
of two by two matrices
over the field of mod2 residue classes. Fix a set of matrix units
$\{e_{ij}\}$ and introduce the coproduct $\cop(e_{ij})\colon
 = e_{ij} \o e_{ij}$. Then we have two  normalized left integrals
$l_j = \sum_i e_{ij}$
for $j=1,2$ neither of which is non-degenerate. The only non-degenerate
left integral is  $l=l_1+l_2$ for which however $\PL(l)=0$. }

As an immediate consequence of the above considerations we have
\begin{scho} \label{ndeg} The following properties for $l\in
\LI$ ($r\in\RI$) are equivalent:
\begin{description}
\item[i)] $l$ ($r$) is non-degenerate;
\item[ii)] $l$ is separating for ${\cal I}^L_{A^{L,R}}$
            ($r$ is separating for $_{A^{L,R}}{\cal I}^R$ );
\item[iii)] $l$ is cyclic for ${\cal I}^L_{A^{L,R}}$
             ($r$ is cyclic for $_{A^{L,R}}{\cal I}^R$).
\end{description}
\end{scho}
 
In a Frobenius
WHA $A$ the group of invertible elements $A^R_\x$ of $A^R$ acts on
the set $\LI_*(A)$ of non-degenerate left integrals transitively
and freely. Similar statement holds for the non-degenerate right
integrals $\RI_*$,
\beq
\LI_*=l\,A^R_\x\ ,\qquad \RI_*=A^L_\x\,r
\eeq
for any $l\in\LI_*$ and $r\in\RI_*$. Similar relation for the dual
integrals shows that there are one-to-one correspondences between
non-degenerate integrals of $A$ and of $\duA$. The Theorem below
selects a distinguished "natural" one-to-one correspondence.
\begin{thm}                 \label{thm: dualint}
Let $A$ be a WHA and $l\in\IL$ be a left integral. If there
exists a $\lambda\in\duA$ such that $\lambda\la l=\1$ then
it is unique, it is a left integral in
$\duA$, and both $l$ and $\lambda$ are
non-degenerate. Moreover $l\la\lambda=\du1$. Such a pair
$(\l,\lambda)$ will be called a dual pair of left integrals.
 
Similarly, elements $r\in\RI_*$ and $\lambda\in\duLI_*$ are in
one-to-one correspondence by either one of the equivalent
relations $\lambda\la r=\1$ or $\lambda\ra r=\du1$.
\end{thm}
 
\Proof
By Lemma \ref{scho: IL}.d) if $l$ is a left integral such
that $\lambda\la l=\1$ then\footnote{Here
we use the standard notations $f_L,f_R\colon
A\to\duA$ defined by $f_L(x):=f\ra x$ and $f_R(x):=x\la f$ for
any $f\in\duA$.}
$l_R\circ\lambda_L=S$. Since $S$ is
invertible, both $l_R$ and $\lambda_L$ are invertible, i.e. $l$
and $\lambda$ are non-degenerate and $\lambda$ is unique. To show
that $\lambda\in\duIL$
\[
\varphi\lambda\la l\ =\ \varphi\la\1\ =\ \duPL(\varphi)\lambda\la
l\ ,\quad \varphi\in\duA
\]
suffices since $l_R$ is a bijection. It remains to show that
$l\la\lambda=\du1$ which eventually justifies the term "dual" left
integral. For $l\in\IL$ and $\lambda\in\duA$ we have
\beanon
\lambda\la l=\1&\Leftrightarrow&
xl\c\bra\lambda,l\cc\ket=x\ \ x\in A\\
&\Leftrightarrow&
\bra\lambda,xl\cc\ket S^{-1}(l\c)=x\ \ x\in A\\
&\Rightarrow&
\bra l\la\lambda,x\ket=\eps(x)\ \ x\in A\\
&\Leftrightarrow&
l\la\lambda=\du1
\eeanon
 
The duality between $\duLI$ and $\RI$ follows from the above
duality between $\duLI$ and $\LI$ by passing from $A$ to $A^{op}$.
The other two twisted versions of the Theorem are not spelled out
explicitely. They can also be obtained by applying the antipode
to the above relations. \qed

Recall that the {\em quasibasis} of a non-degenerate functional
$f$ on $A$ is an element $\sum_i a_i\o b_i\in A\o A$ such that
(cf. \cite{Watatani})
\beq
\sum_i f(xa_i)b_i=x=\sum_i a_i f(b_i x)\,,\quad x\in A\ .
\eeq
(If $K$ is a field then this just means that $\{b_i\}$ is a
$K$-basis of $A$ and $\{a_i\}$ is its dual basis w.r.t. $f$.)
In other words $\sum_i a_i\o b_i$ is simply the expression $\sum_i
f_R^{-1}(\beta^i)\o b_i$ of the inverse of $f_R\colon A \to \duA$
as an element of $A\o A$.
The index of $f$ is then defined by $\Ind f:=\sum_i a_ib_i$ which
belongs to $\Center A$.
Now let $(l,\lambda)$ be a dual pair of left integrals.
Then the quasibasis of $\lambda$ is $l\cc\o S^{-1}(l\c)$ and
\beq
\Ind\lambda=S^{-1}\circ\PL(l)\ \in Z^R\ .
\eeq
In particular a non-degenerate left integral $l$ is normalized if
and only if its dual has index $\1$.

\subsection{2-sided non-degenerate integrals}
The space of 2-sided integrals $\I(A):=\LI(A)\cap\RI(A)$ in a weak
Hopf algebra $A$ is a possibly zero subalgebra of $A$. The
assumption $\I(A)\neq 0$ is independent of the
assumption $\LI_*(A)\neq\emptyset$ since already Hopf algebras
provide examples \cite{Van Daele} for $\LI_*(A)\neq 0$ and $\I(A)=0$.
In this Subsection we will make the
stronger assumption $\I_*(A):=\LI_*(A)\cap\I(A)\neq\emptyset$ and
study some of the consequences. The main result will be finding a
criterion for a WHA to be a symmetric algebra.
 
At first we observe that if a non-degenerate 2-sided integral
$j$ exists
then the subspace of 2-sided integrals is obtained from $j$ by
the action of the central subalgebra $Z^R=A^R\cap\CA$,
\beq
\I=j\,Z^R\,,\quad\I_*=j\, Z^R_\x\qquad\mbox{for any}\quad
j\in\I_*\ .
\eeq
As a matter of fact if $i\in\I$ then
$i$ is a left integral therefore there exists an $x^R\in\A^R$ such
that $i=jx^R$. Thus for all $y\in A$ we have
$jx^R\PR(y)=jx^Ry=j\PR(x^Ry)$. Since $j$ is separating for the
right $A^R$-action, $x^R\PR(y)=\PR(x^Ry)$. Therefore
\[
x^RS(y)=x^RS(y\c)y\cc S(y\ccc)=S(y\c)x^Ry\cc S(y\ccc)=S(\1\c
y)x^R\1\cc=S(y)x^R
\]
hence $x^R$ is central.
 
Next we recall some facts about "modular automorphisms". Let $A$
be a finite dimensional Frobenius algebra over a field $K$ and let
$f\colon A\to K$
be a non-degenerate functional. Then the {\em modular
automorphism}
of $f$ is defined to be the unique $\theta_f\in\Aut A$ such that
\beq
f(xy)\ =\ f(y\theta_f(x))\ ,\qquad x,y,\in A\ .
\eeq
It is worth to give two other equivalent definitions of $\theta$:
\beq
f\ra x\ =\ \theta_f(x)\la f\ ,\qquad x\in A
\eeq
or simply
\beq
\theta_f\ =\ f_R^{-1}\circ f_L\ .
\eeq
Since any two non-degenerate functionals $f$ and $g$ are related
by $g=x\la f$, with $x\in A_\x$, the equivalence class
$\theta_A:=[\theta_f]$ of $\theta_f$ modulo inner automorphisms is
independent of the choice of $f$. If $A$ is a WHA which is
Frobenius then one may ask the question whether $\theta_A=[S^2]$.
\begin{defi}
A non-degenerate functional $f\colon A\to K$ over a WHA $A$ is
called a q-trace if $\theta_f=S^2$.
\end{defi}
In the term "q-trace" the letter "q" has no individual meaning.
One may as well read it as "skew trace" although we do not denie
that our motivation came from the theory of q-deformed Hopf
algebras.
\begin{lem}                          \label{lem: Sinv&qtr}
In a WHA $A$ let $l$ be a non-degenerate left integral. Then
$S(l)=l$ if and only if its dual left integral $\lambda$ is a
q-trace.
\end{lem}
\Proof $\theta_{\lambda}=S^2$ is equivalent to that the quasibasis
of $\lambda$ satisfies
\beq
l\cc\o S^{-1}(l\c)\ =\ S(l\c)\o l\cc\ .
\eeq
Applying $S$ to the second tensor factor we obtain
$\cop(l)=\cop(S(l))$ which yields $l=S(l)$ by the existence of
a counit. \qed
 
\begin{lem}                         \label{lem: S-inv}
If non-degenerate 2-sided integrals exist then all 2-sided
integrals $i\in\I(A)$ are $S$-invariant, $S(i)=i$.
\end{lem}
\Proof If we can show only that the non-degenerate 2-sided
integrals are $S$-invariant then we are ready since
$j=S(j)\in\I_*$ implies $S(jz^R)=S(z^R)j=z^Rj=jz^R$ for all
$z^R\in Z^R$.
 
So let $j\in\I_*$. Then $S(j)\in\I_*$ thus there exists an
invertible $z\in Z^R$  such that $S(j)=jz$.
Let $\lambda$ be the dual of $j$ as a left integral. Then
for arbitrary $x\in A$ and for $z^L=S^{-1}(z^{-1})$
\beanon
z^LS(x)&=&z^L(\lambda\ra x)\la j=(\lambda\ra x)\la z^Lj=
        (\lambda\ra x)\la S^{-1}(j)\\
S^2(x)z^{-1}&=&j\ra\duS^{-1}(\lambda\ra x)=j_L\circ
\duS^{-1}\circ\lambda_L(x)=\lambda_R^{-1}\circ\lambda_L(x)=\\
&=&\theta_\lambda(x)
\eeanon
Therefore $z^{-1}=\theta_\lambda(\1)=\1$ and $j$ is $S$-invariant.
\qed
 
\begin{thm} The WHA $A$ over $K$ is a symmetric algebra if and
only if it has non-degenerate
2-sided integrals and the square of the antipode is an inner automorphism.
\end{thm}
 
\Proof Let $A$ be a symmetric WHA, and $\tau\in {\hat A}$ be a
non-degenerate trace. Then there exists a unique $i\in A$ such that
$i\la \tau =\du1 =\tau\ra i$. We claim that $i$ is a 2-sided
integral. As a matter of fact
\beanon xi\la \tau &=&x\la\du1\ =\ \PL(x)\la\du1\ =\
\PL(x)i\la \tau  \\
\tau \ra ix&=&\du1\ra x\ =\ \du1\ra\PR(x)\ =\
\tau \ra i\PR(x)
\eeanon
so by non-degeneracy of $\tau$, $i\in\I$. This integral $i$ is also
non-degenerate: For any $x^R\in A^R$ one has $ix^R\la \tau =
i\la\tau\ra x^R={\hat \1} \ra x^R$,
hence $i$ is separating for ${\cal I}^L_{A^R}$ so non-degenerate by
Scholium \ref{ndeg}.
 
The innerness of $S^2$ in a symmetric algebra follows if we can
construct a non-degenerate functional on $A$ the modular automorphism
of which is $S^2$. By Lemma \ref{lem: S-inv} $i$ is $S$-invariant so
by Lemma \ref{lem: Sinv&qtr} $\chi$, the dual left integral to
$i$, is such a non-degenerate q-trace.
 
Conversely, let $S^2=\Ad_g$ with some $g\in A_\x$ and $i\in\I_*$.
Denoting the dual left integral of $i$ by $\chi$ again, $g^{-1}\la
\chi$ is a non-degenerate trace.
\qed
 
We close this Subsection with a result arising from assuming the
existence of non-degenerate 2-sided integrals in both $A$ and
$\duA$. Although the arising structure is reminiscent to that of
the "distinguished grouplike element" in Hopf algebra theory it is
not a generalization of that.
\begin{prop}
Let $A$ be a WHA and assume that both $\I_*(A)$ and $\I_*(\duA)$
are non-empty. Then $S^4$ is inner and the square of $\theta_A$ is
the identity in $\Out A$. Moreover and more explicitely, for $\hat
h\in\I_*(\duA)$ there exist invertible elements $a_L\in A^L$ and
$\alpha_L\in\duA^L$ such that, with the notations $a_R=S(a_L)$ and
$\alpha_R=\duS(\alpha_L)$, we have
\bea
\Ad_{a_La_R^{-1}}&=&S^4\label{eq: S^4}\\
\Ad_{a_La_R}&=&\theta_{\hat h}^2\label{eq: theta^2}\\
a_La_R^{-1}\la \psi\ra a_La_R^{-1}&=&\alpha_L\alpha_R\psi
\alpha_R^{-1}\alpha_L^{-1}\,,\quad\psi\in\duA\
\label{eq: a-alfa}
. \eea
\end{prop}
\Proof Choose $h\in \I_*(A)$ and $\hat h\in\I_*(\duA)$ and let
$\lambda$ be the dual of $h$ and $l$ be that of $\hat h$, as left
integrals. Define
\beq
a_L=\hat h\la h\ ,\qquad \alpha_L=h\la\hat h\ .
\eeq
Then
\beq
\du1\ra a_L=\bra\du1\c\hat h,h\ket\du1\cc=\bra\hat h\c,h\ket
\hat h\cc\duS(\hat h\ccc)=\duPL(\hat h\ra h)=\duS(\hat h\ra h)
=\alpha_L
\eeq
and introducing $a_R$ and $\alpha_R$ as above
\beq
la_L=l\ra \alpha_R=l\ra(\hat h\ra h)=h(l\ra\hat h)=h(\duS^2(\hat
h)\la l)=h
\eeq
where q-trace property of $l$ and $\duS$-invariance of $\hat h$
have been used. Similarly,
\bea
\alpha_R\la l&=\ h\ =&la_R\\
a_R\la\lambda&=\ \hat h\ =&\lambda\alpha_R\\
\lambda\ra a_R&=\ \hat h\ =&\lambda\alpha_L .
\eea
Non-degeneracy of $\hat h$ and $h$ now imply invertibility of
$a_L$, $a_R$, $\alpha_L$, and $\alpha_R$. Hence Eq(\ref{eq: a-alfa})
readily follows.
 
We can now compute the modular automorphism of $\hat h$
using the information $\theta_\lambda=S^2$.
\beq
\hat h\ra x=a_R\la\lambda\ra x=
a_R\theta_\lambda(x)\la\lambda\quad\Rightarrow\quad
\theta_{\hat h}=\Ad_{a_R}\circ S^2
\eeq
Computing $\duS(\hat h\ra x)$ in two different ways
\beanon
\duS(\hat h\ra x)&=&S^{-1}(x)\la\hat h=\hat h\ra\theta_{\hat
h}^{-1}(S^{-1}(x))\\
&=&\duS(\theta_{\hat h}(x)\la \hat h)=\hat h\ra
S^{-1}(\theta_{\hat h}(x))
\eeanon
yields
\beq
S^{-1}\circ\theta_{\hat h}=\theta_{\hat h}^{-1}\circ S^{-1}
\eeq
and finally
\beq
\Ad_{a_L}\circ S^{-2}\ =\ \theta_{\hat h}\ =\ \Ad_{a_R}\circ S^2
\eeq
from which (\ref{eq: S^4}) and (\ref{eq: theta^2}) follow
immediately. \qed
 
\subsection{Haar integrals}
Since finite dimensional weak Hopf algebras do not go beyond the
"compact" and "discrete" case, the following very conservative
definition of Haar measure will suffice.
\begin{defi}  \label{def: Haar}
An element $h$ of a WHA $A$ is called a Haar integral in $A$ or
Haar measure on $\duA$ if $h$ is a normalized 2-sided integral,
i.e. $h\in \I(A)$ and $\PL(h)=\PR(h)=\1$.
\end{defi}
Obviously, if Haar integral exists then it is a unique
$S$-invariant idempotent.
As a matter of fact let $h$ and $h'$ be Haar integrals. Then
$h'=\PL(h)h'=hh'=h\PR(h')=h$. In particular $h^2=h$.
$S$-invariance follows from uniqueness since $S(h)$ is a always a
Haar integral if $h$ is.
 
In finding criteria for the existence of Haar measure in $A$ an
important role will be played by a special element $\chi\in\duA$
the definition of which was inspired by similar computations in
Hopf algebra theory \cite{Van Daele}:
\beq     \label{eq: khi}
\chi\ :=\ \sum_i\ \beta^i\ra S^{-2}(b_i)\ \equiv\ \hat L'(\du1)\,,
\eeq
where $\{b_i\}$ and $\{\beta^i\}$ are dual bases of $A$ and
$\duA$, respectively, and $\hat L'\colon\duA\to\duA$ is given by
$\hat L'(\psi):=\sum_i \beta^i\psi\ra S^{-2}(b_i)$ .  Note that
$\hat L'$ is the "cop" version of the dual analogue $\hat L$
of the projection (\ref{eq: L}) onto the space of left integrals.
Hence $\chi$ is a left integral in $(\duA)_{cop}$ and therefore in
$\duA$.
 
As we will see below if $\chi$ is non-degenerate and a q-trace then
its dual left integral will automatically be the Haar measure.
In order to see that it is a q-trace let $\Tr_A$ be the standard
trace on $\End_K A$ and introduce the notation $Q_-(x)y:=yx$.
Then for $x\in A$ we have
\beanon
\chi(x)&=&\sum_i\bra\beta^i,S^{-2}(b_i)x\ket=\Tr_A\ Q_-(x)\circ
S^{-2}\\
\chi(xy)&=&\Tr_A\ Q_-(y)\circ Q_-(x)\circ S^{-2}=
\Tr_A\ Q_-(y)\circ S^{-2}\circ Q_-(S^2(x))=\\
&=&\chi(y S^2(x)).
\eeanon
The next Lemma will be crucial in deciding whether $\chi$ is
non-degenerate.
\begin{lem}                             \label{lem: khi}
Let $l$ be a left integral in a WHA $A$ and let $\chi\in\duA$ be
the q-trace left integral defined in Eq(\ref{eq: khi}). Then
\beq                              \label{eq: l->khi}
l\la\chi\ =\ \duS^2(\du1\ra l)\ .
\eeq
\end{lem}
\Proof  Using the q-trace property of $\chi$ and then (\ref{piba})
\beanon
l\la\chi
&=&\sum_i\beta^i\ra S^{-2}(b_il)=(\du1\ra\1\c)\ra S^{-2}(\1\cc
l)=\\
&=&\bra\du1\c,S^{-1}(\1\cc)\1\c\ket\du1\cc\ra S^{-2}(l)=\du1\ra
S^{-2}(l)=\duS^{2}(\du1\ra l)
\eeanon
\qed
\begin{prop}
Let $A$ be a weak Hopf algebra over a field $K$ and let $\chi$ be
given by (\ref{eq: khi}).
\begin{description}
\item[i)] The Haar integral $h\in A$ exists if and only if $\chi$
is non-degenerate, in which case $(h,\chi)$ is a dual pair of left
integrals. In particular Haar integrals are non-degenerate.
\item[ii)] A left integral $l\in\IL $ is a Haar integral if
and only if $\PR(l)\ =\1$.
\end{description}
\end{prop}
The characterization of Haar measures under (ii) is so simple that
it
could be well used as a definition of Haar measure. Notice that in
that case the formal difference between the notions of normalized
left integral and Haar measure  were so tiny (change $\PL$ for
$\PR$) that it would smear out the big conceptual difference:
The existence of normalized left integrals is equivalent to
semisimplicity while the existence of Haar measures is much
stronger.
 
\Proof {\bf ii)} Assume $l\in\IL$ satisfies $\PR(l)=\1$. Then by
Lemma
\ref{lem: khi} $l\la\chi=\du1$. Therefore the duality Theorem
(Thm \ref{thm: dualint}) implies that $(l,\chi)$ is a dual pair of
non-degenerate left integrals. Since $\chi$ is a q-trace, Lemma
\ref{lem: Sinv&qtr} shows that $l$ is an $S$ invariant
non-degenerate left integral. Furthermore
$\PL(l)=\PL(S(l))=S\circ\PR(l)=\1$. Thus $l$ is a Haar integral.
Now assume $h$ is a Haar integral. Then obviously $h$ is a left
integral satisfying $\PR(h)=\1$.
 
{\bf i)} The "only if" part follows from the proof of (ii). Assume
$\chi$ is non-degenerate and let $h$ be its dual left integral.
Then by Lemma \ref{lem: Sinv&qtr} $h$ is 2-sided and by Lemma
\ref{lem: khi} it is normalized. \qed
 
However simple, the criteria of the above Proposition are very
difficult
to verify in concrete situations. So it is worth looking for other
criteria even if they are not applicable in full generality.
 
\begin{thm}         \label{thm: Haar}
Let $A$ be a WHA over an algebraically closed field $K$. Then a
necessary and sufficient condition for the existence of Haar
measure $h\in A$ is that  $A$ is semisimple and there exists a $g\in A_\x$
such that $gxg^{-1}=S^2(x)$ for $x\in A$
and $\tr\,D_r(g^{-1})\neq 0$ for all irreducible representation
$D_r$ of $A$\,.
\end{thm}
The assumption on $K$ is used only to ensure that $A$ is split
semisimple, $A=\oplus_r M_{n_r}(K)$, once knowing that it is
semisimple. In particular there will be a $K$-basis
$\{e_r^{\alpha\beta}\}$ for $A$ obeying matrix unit relations.
 
\Proof {\em Sufficiency}\,: Let $\tau\colon A\to K$ be the trace
with trace vector $\tau_r=\tr D_r(g^{-1})$. Then $\tau$ is
non-degenerate and has as quasibasis the element
\[
\sum_i x_i\o y_i\ :=\ \sum_r\ {1\over
\tau_r}\sum_{\alpha,\beta=1}^{n_r}\,e_r^{\alpha\beta}\o
e_r^{\beta\alpha}\ .
\]
Notice that $\sum_i x_ig^{-1}y_i=\1$. Now we
define $\chi':=g\la\tau$ and claim that $\chi'$ coincides with the
$\chi$ of Eq(\ref{eq: khi}). As a matter of fact
\beanon
\chi'(x)&=&\tau(gx)=\sum_i\tau(x_ig^{-1}y_igx)=\sum\tau(x_i
S^{-2}(y_i)x)=\\
&=&\sum_i\bra\beta_i,S^{-2}(b_i)x\ket=\chi(x)
\eeanon
where we used the fact that the dual of the basis $b_i=y_i$ is
$\beta_i=\tau\ra x_i$. Since $\chi'$ was non-degenerate by
construction, we conclude that {\em the} $\chi$ of Eq(\ref{eq:
khi}) is non-degenerate and therefore its dual left integral $l$
has $\PR(l)=\1$ by Eq(\ref{eq: l->khi}). Therefore $l$ is a Haar
measure.
 
{\em Necessity}\,: If $h\in A$ is a Haar measure then $A$ is
semisimple
by Theorem \ref{thm: semi}. Therefore $A$ is a symmetric algebra
and $\theta_A=\id$. This means that $\theta_\psi$ is inner for all
non-degenerate functional $\psi$. In particular $\theta_\chi=S^2$
is inner where $\chi$ is the dual left integral of $h$. Choose a
$g\in A_{\x}$ implementing $S^2$ and construct the non-degenerate
trace $\tau:=g^{-1}\la\chi$.
\[
\tau(x)=\chi(xg^{-1})=\Tr_A\ Q_-(xg^{-1})\circ S^{-2}=
\Tr_A\ Q_+(g^{-1})\circ Q_-(x)
\]
where $Q_+(x)y:=xy$ is left multiplication on $A$. Choosing a
matrix unit basis to evaluate the trace we obtain
\[
\tau(x)=\sum_r\ \tr D_r(g^{-1})\ \tr D_r(x)
\]
and by non-degeneracy of $\tau$ all components $\tr D_r(g^{-1})$
of the trace vector are non-vanishing. \qed

\sec{$C^*$-Weak Hopf Algebras}
 
In this Section we introduce the $C^*$-structure in WHAs
which is inevitable if WHAs are to be used as symmetries of
inclusions of von Neumann algebras, in particular in quantum field
theory. Utilizing the results of Sections 2 and 3 we establish the
existence of two canonical elements in any $C^*$-WHA, the Haar
measure $h$ and the canonical grouplike element $g$. While the
Haar measure is well known for $C^*$-Hopf algebras, the
canonical grouplike element cannot be recognized in finite
dimensional Hopf
algebras because it is always equal to $\1$. This is related to
involutivity of the antipode in finite dimensional $C^*$-Hopf
algebras \cite{W}. The very fact that $C^*$-WHAs can have
non-involutive antipodes provides the sufficient flexibility
for the emergence of non-integer dimensions.

\subsection{First consequences of the $C^*$-structure}
\begin{defi}
We define a $^*$-WHA as a WHA $(A,\1,\cop,\eps,S)$ over the
complex
numbers $\C$ together with an antilinear involution $^*$ such that
\begin{description}
\item[i)]\ $(A,\ ^*)$  is a $^*$-algebra,
\item[ii)]\ $\cop$ is a $^*$-algebra map, i.e.
      $(x^*)\c\o(x^*)\cc=(x\c)^*\o(x\cc)^*$ for all $x\in A$.
\end{description}
\end{defi}
By uniqueness of the unit, counit, and the antipode (see Lemma
\ref{lem: uniq S}) we have the following additional relations in a
$^*$-WHA.
\beq   \label{eq: 66}
\1^*=\1,\quad\eps(x^*)=\overline{\eps(x)},\quad
S(x^*)^*=S^{-1}(x).
\eeq
Now it is easy to check that the projections $\PL$ and $\PR$
satisfy
\beq
\PL(x)^*=\PL(S(x)^*),\qquad \PR(x)^*=\PR(S(x)^*)
\eeq
therefore $A^L$ and $A^R$ are $^*$-subalgebras of $A$.
As an elementary exercise we obtain selfduality of the $^*$-WHA:
\begin{scho}
Let $A$ be a $^*$-WHA and define a star operation on its dual as
follows
\beq
\bra\varphi^*,x\ket\ =:\ \overline{\bra\varphi, S(x)^*\ket}\ .
\eeq
Then $\duA$ with this star operation becomes a $^*$-WHA.
\end{scho}
 
For a $^*$-WHA $A$ the canonical isomorphisms $\kappa_A^L\colon
A^L\to\duA^R$ and $\kappa_A^R\colon A^R\to\duA^L$ of Lemma
\ref{kappa} become $^*$-algebra isomorphisms.
 
We omit the discussion of further properties of $^*$-WHAs and turn
to the most important case of $C^*$-WHAs.
\begin{defi}
A $^*$-WHA $A$ possessing a faithful $^*$-representation is called
a $C^*$-weak Hopf algebra, or $C^*$-WHA for short.
\end{defi}
Being a finite dimensional $C^*$-algebra any $C^*$-WHA can be
uniquely characterized, as an algebra, by the dimensions
$n_r\in\N$ of its blocks where $r$ is running over the finite set
$\Sec A$ of equivalence classes of irreducible representations
(i.e. the sectors) of $A$.
\beq
A\ \cong\ \bigoplus_{r\in\Sec A}\ M_{n_r}\ ,\qquad
M_{n_r}=\mbox{Mat}(n_r,\C)\ .
\eeq
$A^L$ and $A^R$ are unital $^*$-subalgebras therefore they are
$C^*$-algebras as well and we have natural numbers $n_a$,
$a\in\Sec A^L$ and $n_b$, $b\in\Sec A^R$ characterizing the type
of $A^L$ and $A^R$, respectively.
\beq
A^c\ \cong\ \bigoplus_{a\in\Sec A^c}\ M_{n_a}\ ,\qquad c=L, R\ .
\eeq
The antiisomorphism $S\colon A^L\to A^R$ establishes a bijection
$a\mapsto\bar a$ of the blocks of $A^L$ to the blocks of $A^R$
such that $n_{\bar a}=n_a$. (We consider $\Sec A^L$, $\Sec A^R$,
and $\Sec A$ as disjoint sets which allows to use one function
$n$.)
 
The following elementary but important Proposition will be the
basic ingredient in proving both the existence of Haar measures
and rigidity of the representation category of $C^*$-WHAs.
\begin{prop}                             \label{pro: g}
Let $A$ be a finite dimensional $C^*$-algebra and $S\colon A\to
A^{op}$ an algebra isomorphism such that $(\,^*\,\circ
S)^2=\id_A$. Then there exists $g\in A_\x$ such that
\begin{description}
\item[i)]$\ g\geq\ 0$
\item[ii)]$\ gxg^{-1}\ =\ S^2(x),\ \ x\in A$
\item[iii)]$\ \tr_r(g)\ =\ \tr_r(g^{-1}),\ \ r\in\Sec A$
\item[iv)]$\ S(g)\ =\ g^{-1}$
\end{description}
where $\tr_r$ denotes trace in the irreducible representation
$D_r$.
An element $g\in A$ satisfying only the first three properties is
already unique.
\end{prop}
\Proof The restriction $S|_{{\rm \scriptstyle Center}\,A}$ is an
algebra automorphism
therefore acts on the minimal central idempotents $e_r$ as
$S(e_r)=e_{\bar r}$ where $r\mapsto \bar r$ is a permutation of
$\Sec A$. Since $e_r^*=e_r$ and $^*\circ S$ is an involution,
$r\mapsto\bar r$ is an involution.
 
Choose matrix units $\{e_r^{\alpha\beta}\}$ for the $C^*$-algebra
$A$ and define the antiautomorphism $S_0\colon A\to A$ by
$S_0(e_r^{\alpha\beta}):=e_{\bar r}^{\beta\alpha}$. Then
$S_0^2=\id_A$ and $^*\circ S_0=S_0\circ\,^*$. Since $S\circ S_0$
is an automorphism of $A$ that acts as the identity on the centre,
there exists $C\in A$ invertible such that $S=\Ad_C\circ S_0$..
It follows that
\beq
^*\circ S(x)=C^{-1*}S_0(x^*)C^*,\quad
(^*\circ S)^2(x)=C^{-1*}S_0(C^{-1})xS_0(C)C^*=x,
\eeq
therefore $S_0(C)C^*$ is central and so is its adjoint
$K:=CS_0(C^*)=S(C^*)C$.
\beq
S^2(x)=CS_0(CS_0(x)C^{-1})C^{-1}=CS_0(C^{-1})xS_0(C)C^{-1},\quad
x\in A
\eeq
hence $T:=CS_0(C^{-1})=CC^*[S_0(C)C^*]^{-1}=CC^*K^{-1*}$
implements $S^2$ and its polar decomposition takes the form
\beq
T=ug'\,,\qquad u=K^{-1*}(K^*K)^{1/2}\,,\quad g'=C(K^*K)^{-1/2}C^*
\ .
\eeq
Using centrality of the unitary part and the computations
$S(T)=S(C)S^2(C^{-1})=$\break $S(C)TC^{-1}T^{-1}=T^{-1}$ and
$S(K)=S_0(K)=C^*S_0(C)=S_0(C)C^*=K^*$ we obtain that $g'$
is positive invertible, implements $S^2$, and satisfies
$S(g')=g'^{-1}$.
These latter three properties, however, do not fix $g$ completely.
If $c$ is positive, central, and satisfies $S(c)=c^{-1}$ then
$g=g'c$ will also satisfy the above three properties. Now defining
\beq
g:=g'c\quad\mbox{where}\quad c=\sum_r\ e_r\left(
{\tr_r(g'^{-1})\over \tr_r(g')}\right)^{1/2}
\eeq
it is easy to verify that $g$ obeys (i--iv)
of the Proposition. If $f\in A$ satisfies only (i),
(ii), and (iii) then $f=gc$
where $c$ is positive invertible, central, and
satisfies $D_r(c)=D_r(c)^{-1}$ for all irrep $D_r$. Hence $c=\1$,
proving uniqueness of $g$. \qed
 
\subsection{The Haar measure and selfduality}
Recall that the Haar measure in a WHA $A$ has been defined in
Definition \ref{def: Haar} as the unique element $h\in A$ making
the integral $\int\varphi:=\bra\varphi,\,h\ket$ of a function
$\varphi\colon A\to\C$ to be a non-degenerate functional invariant
under left and right translations and normalized according to
$\int\varphi^L=\dueps(\varphi^L)$ for $\varphi^L\in\duA^L$. The
sufficient conditions for its existence given by Theorem \ref{thm:
Haar} will be used here to prove the next Theorem.
\begin{thm}             \label{thm: C^* Haar}
In a $C^*$-WHA $A$ Haar measure $h\in A$ exists. It is
selfadjoint, $h^*=h$, and such that
\beq
(\varphi,\psi)\ :=\ \bra\varphi^*\psi,\ h\ket\ ,
\qquad \varphi,\psi\in\duA\ ,
\eeq
is a scalar product on $\duA$ making $\duA$ a Hilbert space and
making the left regular module $_\duA\duA$ a faithful
$^*$-representation of the $^*$-WHA $\duA$. Thus $\duA$ is a
$C^*$-WHA, too.
\end{thm}
\Proof $A$ being a finite dimensional $C^*$-algebra is semisimple.
By Proposition \ref{pro: g} there exists a $g$ implementing
$S^2$. This $g$ was shown to be positive and invertible, hence
$\tr D_r(g^{-1})>0$ for all $r\in\Sec A$. Therefore all the
conditions of Theorem \ref{thm: Haar} are satisfied and Haar
measure $h$ exists.
 
Since $h$ is non-degenerate, $(\ ,\ )$ is a non-degenerate
sesquilinear form on $\duA$. So it remains to show positivity.
By the equality
\beq
(\psi,\psi)=\bra\psi^*\psi,\,h\ket=\overline{\bra\psi,S(h\c)^*\ket}
\ \bra\psi,h\cc\ket
\eeq
positivity of $(\ ,\ )$ follows if we can show that
$(S\o\id)\circ\cop(h)$ belongs to the positive cone
\beq  \label{Pcone}
{\cal P}\ =\ \left\{\sum_k a_k^*\o a_k\ |\ a_k\in A\ \right\}\
\subset A\o A\
. \eeq
Therefore the next Lemma will complete the proof. \qed
\begin{lem}                    \label{lem: coph}
Choose matrix units $\{e_q^{\alpha\beta}\}$ for $A$ and let
$g$ denote the element determined in Proposition \ref{pro: g}. If
furthermore $\sum_i x_i\o y_i$ is the quasibasis of the trace
$\tau\colon A\to\C$ with trace vector $\tau_q=\tr_q(g^{-1})$
then
\bea                \label{eq: cophR}
S(h\c)\o h\cc&=&\sum_i x_i\o g^{-1}y_i\ =\
\sum_{q\in\Sec A}\ {1\over\tau_q}\sum_{\alpha\beta}
e_q^{\alpha\beta}g^{-1/2}\o g^{-1/2}e_q^{\beta\alpha}\\
h\c\o S(h\cc)&=&\sum_i x_ig\o y_i\ =\ \sum_{q\in\Sec A}
\ {1\over\tau_q}\sum_{\alpha\beta}  \label{eq: cophL}
e_q^{\alpha\beta}g^{1/2}\o g^{1/2}e_q^{\beta\alpha}
\eea
\end{lem}
\Proof The quasibasis of $\chi=g\la\tau$ is $\sum x_ig^{-1}\o y_i$
and since $\chi$ is the dual left integral of $h$,
this quasibasis is equal to $h\cc\o S^{-1}(h\c)$. This implies the
first row. By property {\bf iii)} of Proposition \ref{pro: g}
$\tau$ is an $S$-invariant trace, therefore its quasibasis
can also be written as $\sum_i y_i\o x_i=\sum_i S^{-1}(x_i)\o
S^{-1}(y_i)$. Thus the second row follows from the first. \qed
 
 From now on $h\in A$ will always denote the Haar measure of $A$
and $\hat h\in\duA$ that of $\duA$.
\begin{lem}
In a $C^*$-WHA $A$ the counit is a positive linear functional,
$\eps(x^*x)\geq 0$, $x\in A$.
\end{lem}
\Proof
\[
\eps(x^*x)=\eps(x^*\1\c)\eps(\1\cc\1_{(2')})\eps(\1_{(1')}x)=
\eps(\PL(x)^*\PL(x))=\bra\hat h,\,\PL(x)^*\PL(x)\ket
\geq 0\,,
\]
where we have used $\hat h|_{A^L}=\eps|_{A^L}$, which follows
from $\bra\hat h,\,x^L\ket=\bra\duPL(\hat h),\,x^L\ket
=\bra\du1,\,x^L\ket$ for all $x^L\in A^L$.
\qed
 
$A$ being semisimple the trivial representation $\V$ decomposes
into irreducibles $V_q$ each of them with multiplicity 1 by
Proposition \ref{prop: End trivi}.
The sectors $q\in\Sec A$ occuring in $\V$ with non-zero
multiplicity will be called {\em vacuum sectors}.
\beq
\V\ \cong\ \bigoplus_{q\in\Vac A}\ V_q\ .
\eeq
By Proposition \ref{prop: End trivi} there is a bijection
$q\mapsto z^L_q$
from the set $\Vac A$ of vacuum sectors to the set of minimal
projections in $Z^L$ such that, with $z^R_q:=S(z^L_q)$, we have
\bea
\triv(z^L_q)\ =&\triv(e_q)&=\ \triv(z^R_q)\\
z^L_q=\PL(e_q)&\quad&\PR(e_q)=z^R_q
\eea
where $e_q$ denotes the minimal central projection in $A$
supporting the irreducible vacuum representation $D_q$.
\begin{lem}
$D_r(h)$ is a 1-dimensional projection for $r\in\Vac A$ and
$D_r(h)=0$ if $r$ is not a vacuum sector. The algebra of 2-sided
integrals is generated by minimal projections $h_q$
\beq
\I(A)\ =\ hAh\ =\ \Span\{h_q\,|\,q\in\Vac A\,\}\,,\qquad h_q=he_q
\ .
\eeq
The non-degenerate 2-sided integrals are presisely the invertible
elements: $\I_*(A)=\I(A)_\x$.
\end{lem}
\Proof If $D_r(h)\neq 0$ then pick up a non-zero vector $v_r$ from
the subspace $D_r(h)V_r$ of the irreducible $A$-module $V_r$ and
define
\beq
T\colon A^L\to V_r\ ,\qquad Tx^L:=D_r(x^L)v_r\ .
\eeq
This map is a non-zero left $A$-module map if we equip $A^L$ with
the structure of the trivial $A$-module $_AA^L$ introduced in
Lemma \ref{lem: trivials}. Indeed,
\beq
D_r(x)Tx^L=D_r(xx^Lh)v_r=D_r(\PL(xx^L)h)v_r=T\PL(xx^L)\ .
\eeq
Therefore $r\in\Vac A$.  This proves that $D_r(h)=0$ for
$r\not\in\Vac A$.
 
Now let the Haar integral act on the trivial left $A$-module
$_A\duA^R$.
\beq
\triv(h)\varphi^R\ =\ h\la \varphi^R\ \in\
\duA^L\cap\duA^R\ \equiv\ \hat Z\ .
\eeq
Thus $\triv(h)\colon\duA^R\to\hat Z$ is a projection, onto.
If $z^L$ is a minimal projection in $Z^L$ then $z^L\la\du1$ is a
minimal projection in $\hat Z$ by Lemma \ref{Z}. Hence
$\triv(z^Lh)$
maps $\duA^R$ onto $(z^L\la\du1)\hat Z\cong\C$. This proves that
$\triv(z^Lh)$, the restriction of which is precisely $D_q(h)$ for
some $q\in\Vac A$, is a 1-dimensional projection.
If $i\in\I$ then by the 2-sided normalization of $h$ one can write
$i=hih$. Conversely, $hxh$ is a 2-sided integral for all $x\in A$.
This proves the remaining assertions. \qed
 
The Haar measure provides conditional expectations
\bea   \label{eq: E^LR}
E^L\colon A\to A^L\,,&\quad&E^L(x)=\hat h\la x\\
E^R\colon A\to A^R\,,&\quad&E^R(x)=x\ra \hat h
\eea
As a matter of fact by Lemma \ref{scho: IL}.c) the image of $E^L$
is in $A^L$
since $\hat h$ is a left integral. $E^L$ is unit preserving since
$\hat h$
is normalized. Finally, $E^L$ is positive since $\hat h$ is
positive and $\cop$ is a $^*$-algebra map.

\subsection{The canonical grouplike element}
In this Subsection we investigate further properties of the
element $g$ of Proposition \ref{pro: g}. We show that it is always
a product of left and right elements, implying its
grouplikeness immediately, and obtain expressions for the
modular automorphisms of the Haar measures of $A$ and $\duA$.
 
\begin{prop}        \label{pro: gg}
In a $C^*$-WHA $A$ there exists a unique $g\in A$ such that
\begin{description}
\item[i)] $g\geq 0$ and invertible,
\item[ii)] $gxg^{-1}=S^2(x)$ for all $x\in A$,
\item[iii)] $h\cc\o h\c=h\c\o gh\cc g$.
\end{description}
\end{prop}
\Proof {\em Existence}\,: Let $g$ be the (unique) element defined
by the conditions of Proposition \ref{pro: g}. As in the
proof of Lemma \ref{lem: coph} let $\tau$ be the
$S$-invariant trace
with trace vector $\tau_q=\tr_q(g)$ and $\sum x_i\o y_i$ be its
quasibasis. Then
\bea
h\cc\o h\c&=&\sum_i x_i\o S(g^{-1}y_i)=\sum_i S^{-1}(y_i)\o
x_ig=\\
&=&\sum_i S^{-1}(y_ig^{-1})\o gx_ig=\sum_igS^{-1}(y_i)\o gx_ig
=\\
&=&\sum_i S(g^{-1}y_i)\o gx_ig\ =\ h\c\o gh\cc g
\eea
{\em Uniqueness}\,: Let $g$ and $g'$ satisfy {\bf i)}, {\bf ii)},
and {\bf iii)}. Then $g'=gc$ with $c$ central, positive, and
invertible. Furthermore, since {\bf iii)} is equivalent to
\beq
\bra \varphi\psi,\,h\ket\ =\ \bra\psi(g\la\varphi\ra g),\,h\ket\,,
\eeq
non-degeneracy of $h$ implies
\beq
g'\la\varphi\ra g'\ =\ g\la\varphi\ra g\ ,\qquad \varphi\in\duA\ .
\eeq
Therefore $c^2\la\varphi\equiv c\la\varphi\ra c=\varphi$ for all
$\varphi\in\duA$. Thus $c^2=\1$ and, by positivity, $c=\1$.
\qed
 
Notice that property {\bf iii)} of Proposition \ref{pro: gg} is
equivalent to that the modular automorphism of the
Haar functional $\varphi\mapsto\varphi(h)$ is expressible in the
form
\beq
\theta_h(\psi)\ =\ g\la\psi\ra g\ ,\qquad\psi\in\duA\ .
\eeq
\begin{defi}
Let $A$ be a $C^*$-weak Hopf algebra. Then the unique element
$g\in A$ determined either by the conditions of Proposition
\ref{pro: g} or by the conditions of Proposition \ref{pro: gg}
is called the {\em canonical grouplike element} of $A$.
\end{defi}
As one may suspect the canonical grouplike element is grouplike
in the sense of
\begin{defi}                               \label{grouplikeness}
An element $x$ of a WHA $A$ is called {\em grouplike} if
\bea
\cop(x)&=&x\1\c\o x\1\cc\ =\ \1\c x\o\1\cc x\label{grouplike1}\\
S(x)x&=&\1\ .\label{grouplike2}
\eea
\end{defi}
We remark that if (\ref{grouplike1}) holds then condition
(\ref{grouplike2}) is equivalent to the assumption that $x$ is
invertible.
One should emphasize that grouplike elements are not always like
group elements if a $^*$-operation is present. Namely we allow
for $x$ not to be unitary. Thus there can be positive grouplike
elements, for example, in a $C^*$-WHA.
 
If $x$ is
an invertible element factorizable as $x_Lx_R^{-1}$ with $x_L\in
A^L$ and $x_R=S(x_L)=S^{-1}(x_L)$ then $x$ is automatically
grouplike. As a matter of fact
$\cop(x)=x_L\1\c\o x_R^{-1}\1\cc=xx_R\1\c\o x_R^{-1}\1\cc=x\1\c\o
x\1\cc$. Now it follows from the next Lemma that the canonical
grouplike element $g$ is grouplike.
\begin{lem}
In a weak $C^*$-Hopf algebra $A$ the elements $h\ra\hat h$ and
$\hat h\la h$ are positive and invertible. The canonical grouplike
element of $A$ can be factorized as
\bea
g\ =\ g_Lg_R^{-1}\ ,&\qquad&\mbox{where}\label{eq: g=}\\
g_L:=(\hat h\la h)^{1/2}\ ,&\qquad& g_R=(h\ra\hat h)^{1/2}
\eea
\end{lem}
\Proof $\hat h\la h=E^L(h)=E^L(h^*h)\geq 0$
and similarly $h\ra\hat h\geq 0$
by positivity of the conditional expectations (\ref{eq: E^LR}).
Invertibility follows from the existence of the dual left integral
$\chi$ since $(h\ra\hat h)\la\chi=\duS(\hat h)(h\la\chi)=\hat h$
can hold for the non-degenerate $\chi$ and $\hat h$ only if
$h\ra\hat h$ is invertible. Thus $\hat h\la h=S(h\ra\hat h)$
is invertible, too.
 
The next point is to observe that the three elements $\hat h\la
h$,
$h\ra\hat h$, and $g$ commute with each other. For $g$ and any one
of the others this follows from the fact that $\hat h\la h$ and
$h\ra\hat h$ are invariant under $S^2$. For the commutativity of
the remaining two notice that one of them belongs to $A^L$ the
other to $A^R$. Now compare the following expressions:
\bea
\hat h&=&(h\ra\hat h)\la\chi\ =\ (h\ra\hat h)g\la\tau\ ,\\
\hat h&=&\duS^{-1}(\hat h)=\tau\ra g^{-1}(\hat h\la h)=g^{-1}(\hat
h\la h)\la\tau\ .
\eea
By non-degeneracy of $\tau$ we obtain
\beanon
(\hat h\la h)g^{-1}&=&(h\ra\hat h)g\\
(\hat h\la h)(h\ra\hat h)^{-1}&=&g^2
\eeanon
and taking the (positive) square root the Lemma is proven. \qed
 
\begin{lem}
The left--right components of the canonical grouplike element $g$
of $A$ and $\hat g$ of $\duA$ obey the following identities.
\bea
\hat g_L=\du1\ra g_L=\du1\ra g_R&\qquad&
g_L=\1\ra\hat g_L=\1\ra\hat g_R\\
\hat g_R=g_R\la\du1=g_L\la\du1&\qquad&
g_R=\hat g_R\la\1=\hat g_L\la\1\\
S(g_L)\ =\ g_R\ =\ S^{-1}(g_L)&\qquad&
\duS(\hat g_L)\ =\ \hat g_R\ =\ \duS^{-1}(\hat g_L)
\eea
\end{lem}
\Proof
Since $g_L\in A^L$ and $g_R\in A^R$, they commute and both of them
are invariant under $S^2=\Ad_g$. So are the $C^*$-algebras
generated by each of them, pointwise. Hence
$S(g_L^{1/2})^*=S^{-1}(g_L^{1/2})=S(g_L^{1/2})$ therefore
$S(g_L)=S(g_L^{1/2})^2\geq 0$. On the other hand
$S(g_L)^2=S(g_L^2)=g_R^2$, therefore $S(g_L)$ is the positive
square root of $g_R^2$, i.e. $S(g_L)=g_R$.
 
Next we want to show that $\du1\ra(\hat h\la
h)=h\la\hat h$. Since both hand sides belong to $\duA^L$, the
identity
\beanon
\bra\du1\ra(\hat h\la h),\,x^R\ket&=&\eps((\hat h\la h)S(x^R))=
\eps(\hat h\la hS(x^R))=\\
&=&\bra\hat h\ra h,\,S(x^R)\ket=\bra h\la\hat h,\,x^R\ket\ ,
\eeanon
valid for $x^R\in A^R$, suffices. Therefore $\du1\ra g_L^2=\hat
g_L^2$, or $\du1\ra g_R^2=\hat g_L^2$. Now use the fact that
$A^R\ni x^R\mapsto (\du1\ra x^R)\in\duA^L$ is a $^*$-algebra
isomorphism. Hence passing to the square roots we obtain $\1\ra
g_R=\hat g_L$. All the remaining identities are simple
consequences of this. \qed

\begin{prop}
Let $A$ be a $C^*$-WHA with dual $\duA$ and let $h\in A$, $\hat
h\in\duA$ be the corresponding Haar measures. Then
\begin{description}
\item[i)]
the modular automorphism of the Haar functional $\hat h$ is
implemented by $g_Lg_R$, i.e. for all $x\in A$ we
have $\theta_{\hat h}(x)=g_Lg_R\,x\,g_R^{-1}g_L^{-1}$;
\item[ii)]
the dual left integral of $h$ can be expressed as $\chi=\hat h\hat
g_R^{-2}$;
\item[iii)]
the $S$-invariant trace functional $\tau=g^{-1}\la\chi$ and
the Haar functional $\hat h$ are related by
\bea
\tau&=&\hat g_L^{-1}\hat h\hat g_R^{-1}\\
\hat h&=&g_Lg_R\la\tau
\eea
\end{description}
\end{prop}
\Proof {\bf i)}: Using identities like $\hat g_L\la
x=g_Rx,\dots$etc, which follow from Scholium \ref{sch: arrow}, one
can easily verify $\hat g\la x\ra \hat g=g_Lg_Rxg_R^{-1}g_L^{-1}$,
for $x\in A$.
 
{\bf ii)}: The identity $\hat h\la h=g_L^2=\1\ra\hat g_R^2$
implies $\1=\hat h\la h\ra\hat g_R^{-2}=\hat h\la hg_L^{-2}$,
hence $hg_L^{-2}=hg_R^{-2}$ is the dual left integral of
$\hat h$. By duality, $\hat h\hat g_R^{-2}$ is the dual left
integral $\chi$ of $h$.
 
{\bf iii)}: $\tau=g^{-1}\la\hat h\hat g_R^{-2}=\hat g_R^{-1}(\hat
h\hat g_R^{-2})\hat g_R=\hat g_L^{-1}\hat h\hat g_R^{-1}$ and
$\tau=g^{-1}\la(g_R^{-2}\la\hat h)=g_L^{-1}g_R^{-1}\la\hat
h$ completes the proof. \qed
 
Cyclicity and separability of the vector $h$ in the right
$A^{L,R}$-module $\LI$ (cf. Scholium \ref{ndeg}) allows us to
introduce $\duA^R$-valued "Radon-Nikodym derivatives" of left
integrals $l$ with respect to the Haar measure. At first note that
$l=\PL(h)l=hl=h\PR(l)=hS^{-1}(\PR(l))$ therefore using Scholium
\ref{sch: arrow} we have
\beq
\bra
\varphi,\,l\ket=\bra\varphi\rho_R,\,h\ket=\bra\rho_L\varphi,\,h\ket
\eeq
where $\rho_R=\PR(l)\la\du1$ and
$\rho_L=S^{-1}(\PR(l))\la\du1=\duS^2(\rho_R)$.
 
\begin{prop}
The bijections $\IL\to\duA^R$ provided by the
left and right Radon-Nikodym derivatives $l\mapsto\rho_L$ and
$l\mapsto\rho_R$, respectively, obey the following properties.
\begin{description}
\item[i)] $l$ is non-degenerate iff $\rho_{R,L}$ is invertible.
\item[ii)] If $l$ is non-degenerate then $l$ is normalized iff
$l^2=l$.
\item[iii)] $l$ is of positive type, i.e.
$\bra\varphi^*\varphi,l\ket\geq 0$ for all $\varphi\in\duA$, iff
$\PR(l)\geq 0$ iff $\rho_R\in\hat g_R^{1/2}\duA^R_+
\hat g_R^{-1/2}$ where $\duA^R_+$ is the cone of positive elements
in $\duA^R$. In this case $\rho_L=\rho_R^*$ and there
exists a $\xi\in\duA$ such that
$\bra\varphi,\,l\ket=\bra\xi^*\varphi\xi,\,h\ket$
for $\varphi\in\duA$.
\item[iv)] Let $\lambda$ be the dual left integral of $l$. Then
the Radon-Nikodym derivatives of $\lambda$ and $l$ are related by
$\PR(l)(\duPR(\lambda)\la\1)=g_R^{-2}$.
\end{description}
\end{prop}
\Proof {\bf i)} follows from cyclicity of $h$ in $\LI_{A^R}$.
{\bf ii)}: $l^2=l$ implies $(\PL(l)-\1)l=0$ and acting with
$\lambda\la$, where $\lambda$ is the dual left integral of $l$,
one obtains $\PL(l)=\1$. The converse implication is trivial.
{\bf iii)}: As in the proof of Theorem \ref{thm: C^* Haar} $l$ is
of positive type iff $S(l\c)\o l\cc$ belongs to the positive cone
(\ref{Pcone}). If it does then $\PR(l)=S(l\c)l\cc\geq 0$. Now
assume $\PR(l)\geq 0$. Then introducing $\xi=\PR(l)^{1/2}\la\du1$
we have
$\PR(l)^{1/2}=\xi\la\1$,
$S^{-1}(\PR(l)^{1/2})=S(\xi\la\1)^*=(\1\ra \duS^{-1}(\xi))^*=
\1\ra\xi^*$ therefore $l=hS^{-1}(\PR(l)^{1/2})\PR(l)^{1/2}=\xi\la
h\ra\xi^*$ proving that $l$ is of positive type. It remained to
reformulate positivity of $\PR(l)$ in terms of $\rho_R$.
Use the fact that the antimultiplicative map
$x^R\mapsto(x^R\la\du1)$ from $A^R$ to $\duA^R$ sends
the $^*$-operation into a new involution,
$x^{R*}\la\du1=(S^{-1}(x^R)\la\du1)^*=(S^{-2}(x^R)\la\du1)^*=
(g_Rx^Rg_R^{-1}\la\du1)^*=\hat g_R(x^R\la\du1)^*\hat g_R^{-1}$.
Therefore the equality $\PR(l)=x^{R*}x^R$ for some $x^R\in A^R$ is
equivalent to the equality $\rho_R=(x^R\la\du1)(x^{R*}\la\du1)=
\hat g_R^{1/2}\eta\eta^*\hat g_R^{-1/2}$ with $\eta=\hat
g_R^{-1/2}(x^R\la\du1)\hat g_R^{1/2}\in\duA^R$.
{\bf iv)} follows by an elementary calculus starting from the
identity $\1=\lambda\la l=\hat h\duPR(\lambda)\la
hS^{-1}(\PR(l))$. \qed
\appendix
\sec{Appendix: The Weak Hopf Algebra $B\o B^{op}$}
 
Let $B$ be a separable algebra over the field
$K$ and let $E\colon B\to K$ be a non-degenerate functional with
index $\1$. These are the data needed for constructing a WHA
structure on the algebra $B\o B^{op}$. For a similar construction
of a WBA see \cite{WBA}.
 
At first choose a basis $\{e_i\}$ of $B$ over $K$ and let
$\{f_i\}$ be its dual basis w.r.t. $E$, i.e.
$E(e_if_j)=\delta_{ij}$. Then
\begin{description}
\item[a)] $\sum_i f_i\o e_i\in B\o B$ is independent of the choice
of $\{e_i\}$;
\item[b)] $\sum_i E(xf_i)e_i\ =\ x\ =\ \sum_i f_iE(e_ix)$, $x\in
B$;
\item[c)] $\sum f_ie_i\ =\ \1$;
\item[d)] $\sum_i xf_i\o e_i\ =\ \sum_i f_i\o e_ix$, $x\in B$;
\item[e)] if $\theta$ denotes the modular automorphism of $E$,
i.e. $E(xy)=E(y\theta(x))$, $x,y\in B$, then
\[
\sum_i f_i\o xe_i\ =\ \sum_i f_i\theta(x)\o e_i\,,\quad x\in B\,;
\]
\item[f)] $\sum_i f_i\o e_i\ =\ \sum_i e_i\o \theta^{-1}(f_i)\ =\
\sum_i\theta(e_i)\o f_i$.
\end{description}
The algebra $B\o B^{op}$ is the $K$-space $B\o B$ with
multiplication $(a\o b)(x\o y):=(ax\o yb)$. Its WHA structure is
defined by
\bea
\cop(x\o y)&=&\sum_i\ (x\o f_i)\ \o\ (e_i\o y)\\
\eps(x\o y)&=&E(xy)\\
S(x\o y)&=&y\o \theta(x)
\eea
The verification of the WHA axioms is left to the reader.
The left and right subalgebras of $B\o B^{op}$ are $B\o\1$ and
$\1\o B$, respectively, because we have
\beq
\PL(x\o y)=xy\o \1\,,\qquad\PR(x\o y)=\1\o y\theta(x)\ .
\eeq
Let $A$ be an arbitrary WHA over $K$. Then $A^LA^R$ is a sub-WHA
with hypercenter $A^L\cap A^R$. Thus $A^LA^R$ decomposes into a
direct sum of WHA-s each summand being isomorphic to a WHA of the
type $B\o B^{op}$.
 
Since $B\o B^{op}$ is separable, by Theorem \ref{thm: semi}, it
must contain a normalized left integral. Indeed,
\beq
l:=\sum_i f_i\o e_i\ \equiv \ S^2(\1\cc)\1\c
\eeq
is such a left integral. What is more, it is non-degenerate.
 
Before looking for Haar integrals some remarks about innerness of
$\theta$ are in order. The quantity $q=\sum_i e_if_i$ always
implements $\theta^{-1}$, i.e. $xq=q\theta(x)$ for $x\in B$, but
it is not necessarily invertible. (For example for
$B=M_2(\Z_2)$ and for any non-degenerate functional $E$ the $q$ is
identically zero.) In fact $q$ is invertible iff the left regular
trace on $B$ is non-degenerate (especially if $K$ is of
characteristic zero). Fortunately one can circumvent this nuisance
by using the existence of a non-degenerate trace $\tr$ on any
separable algebra $B$ (see \cite{Curtis-Reiner new}). Then the
Radon-Nykodim derivative $\gamma$ of $E$ w.r.t. $\tr$ provides an
invertible element implementing $\theta$,
\beq
E(x)=\tr(x\gamma)\ ,\qquad \theta(x)=\gamma x\gamma^{-1}\ ,\quad
x\in B\ .
\eeq
This proves that $\theta$ is inner and therefore so is the square
of the antipode, $S^2=\theta\o\theta$.
 
Omitting the details we can now formulate the condition for the
existence of the Haar measures $h$ and $\hat h$ as follows.
Haar measure in $B\o B^{op}$ exists iff $\sum_i f_i\gamma^2
e_i$ is invertible and Haar measure in $\widehat{B\o B^{op}}$ exists iff
$E(\1_B)\neq 0$.


\begin{thebibliography}{XXXXXX}
\addcontentsline{toc}{section}{\protect\numberline{}{References}}
\renewcommand{\b}{\bibitem}
\renewcommand{\baselinestretch}{.3}
\small
\medskip
\b{Abe} E. Abe, {\em Hopf Algebras}, Cambridge University Press,
Cambridge, 1980
\b{BSz} G. B\"ohm, K. Szlach\'anyi, A Coassociative
$C^*$-Quantum Group with Nonintegral Dimensions, {\em Lett. Math.
Phys.} {\bf 35}, 437 (1996)
\b{BNSzII} G. B\"ohm, F. Nill, K. Szlach\'anyi, Weak Hopf
Algebras II: Representation Theory, Dimensions, and Markov
Traces, in preparation
\b{Curtis-Reiner}C. W. Curtis, I. Reiner, {\em Representation
Theory of Finite Groups and Associative Algebras}, John Wiley \&
Sons, Inc., 1962
\b{Curtis-Reiner new} C. W. Curtis, I. Reiner, {\em Methods of
Representation Theory}, John Wiley \& Sons
\b{Drinfeld} V.G. Drinfeld, Quasi-Hopf Algebras, {\em Leningrad
Math. J.} {\bf 1}, 1419-1457 (1990)
\b{FGV} J. Fuchs, A. Ganchev, P. Vecserny\'es,
Rational Hopf algebras: Polynomial equations, gauge fixing,
and low dimensional examples, {\em Int. J. Modern Physics}
{\bf A 10}, 3431 (1995)
\b{Hayashi} T. Hayashi, Quantum group symmetry of partition
functions of IRF models and its application to Jones' index
theory, {\em Commun. Math. Phys.} {\bf 157}, 331-345 (1993)
\b{LS} R. G. Larson, M. Sweedler, An associative orthogonal
bilinear form for Hopf algebras, {\em Amer. J. Math.} {\bf 91},
75-93 (1969)
\b{Longo} R. Longo, A duality for Hopf algebras and for
subfactors I, {\em Commun. Math. Phys.} {\bf 159}, 133-150 (1994)
\b{MS} G. Mack and V. Schomerus, Quasi Hopf quantum symmetry in
quantum theory, {\em Nucl. Phys.} {\bf B370}, 185 (1992)
\b{McLane} S. MacLane, {\em Categories for the Working
Mathematician}, Springer 1971
\b{Montgomery} S. Montgomery, {\em Hopf Algebras and their Actions
on Rings}, CBMS Series of the Amer. Math. Soc., No. 82, 1993
\b{WBA} F. Nill, Axioms for Weak Bialgebras, q-alg/9805104
\b{NSzW} F. Nill, K. Szlach\'anyi, H.-W. Wiesbrock, Weak Hopf
Algebras and Reducible Jones Inclusions of Depth 2. I. From
crossed products to Jones towers, q-alg/9805nnn
\b{Ocneanu 1} A.Ocneanu, Quantized groups, string algebras,
and Galois theory for algebras, in {\em Operator Algebras and
Applications}, Vol. 2, eds.: D.E. Evans {\em et al}., London Math.
Soc. Lect. Notes {\bf 135}, Cambridge 1988
\b{Ocneanu 2} A. Ocneanu, Quantum Cohomology, Quantum Groupoids,
and Subfactors, unpublished talk presented at the First Caribic
School of Mathematics and Theoretical Physics, Guadeloupe 1993
\b{Pierce} R. S. Pierce, {\em Associative Algebras},
(Graduate Texts in Matematics 88) Springer-Verlag 1982
\b{Sweedler} M.E. Sweedler, {\em Hopf algebras}, Benjamin  1969
\b{Sz} K. Szlach\'anyi, Weak Hopf Algebras, in
{\em Operator Algebras and Quantum Field Theory}, eds.: S.
Doplicher, R. Longo, J. E. Roberts, and L. Zsid\'o, International
Press (1996)
\b{Van Daele} A. Van Daele, The Haar measure on finite
quantum groups, to appear in {\em Proc. Amer. Math. Soc.}
\b{V} P. Vecserny\'es, On the Quantum Symmetry of the Chiral Ising
Model, {\em Nucl. Phys.} {\bf B 415 [FS]}, 557-588 (1994)
\b{Watatani} Y. Watatani, Index for $C^*$-subalgebras, {\em
Memoirs of the AMS} No. 424 (1990)
\b{W} S.L. Woronovicz, Compact Matrix Pseudogroups, {\em Commun.
Math. Phys.} {\bf 111}, 613-665 (1987)
\b{Yamanouchi} T. Yamanouchi, Duality for Generalized Kac Algebras
and a Characterization of Finite Groupoid Algebras, {\em Journal
of Algebra} {\bf 163}, 9-50 (1994)
 
\end{thebibliography}
\end{document}